\documentclass[final]{siamltex}
\usepackage{amsmath,graphics,rotating,amssymb,psfrag}
\usepackage[applemac]{inputenc}
 
\numberwithin{remark}{section} 
\newcommand \be {\begin{equation}} 
\newcommand \ee {\end{equation}} 
\newcommand \eps \epsilon
\newcommand \ah {\widehat a} 
\newcommand \uh {\widehat u} 
\newcommand{\RR}{\mathbb{R}}
\newcommand{\BigO}[1]{\mathcal{O}({#1})}

\newcommand{\Lip}{\operatorname{Lip}}
\newcommand{\TV}{ \operatorname{TV}}
\newcommand{\SV}{ \operatorname{SV}}
\newcommand{\mysigma}{\sigma}
\newcommand{\phib}{\varphi^{\flat}}
\newcommand{\phibo}{\varphi^{\flat}_0}
\newcommand{\phibSqr}{\phibo \circ \phib}
\newcommand{\phis}{\varphi^{\sharp}}
\newcommand	\Rp		{{R^\uparrow_+}}
\newcommand	\Cp		{C^{\downarrow}_+}
\newcommand	\CDown	{{C^\downarrow}}
\newcommand	\Cpm 	{{C^\downarrow_\pm}}
\newcommand	\Npm		{{N^\downarrow_\pm}}
\newcommand	\Rm		{{R^\downarrow_-}}
\newcommand	\Cm		{{C^\uparrow_-}}
\newcommand	\CUp		{{C^\uparrow}}
\newcommand	\Cmp		{{C^\uparrow_\mp}}
\newcommand	\Nmp		{{N^\uparrow_\mp}}

\title{Diminishing functionals for nonclassical entropy solutions selected by kinetic relations}

\author{
M{\scriptsize arc} L{\scriptsize aforest}\thanks{D\'epartement de math\'ematiques et g\'enie industriel,
	\'Ecole Polytechnique de Montr\'eal,
	Montr\'eal, Qu\'ebec, Canada, H3C 3A7. E-mail: {\tt Marc.Laforest@polymtl.ca}.}
        \and P{\scriptsize hilippe} G. L{\scriptsize e}F{\scriptsize loch}\thanks{Laboratoire Jacques-Louis Lions
        \& Centre National de la Recherche Scientifique (CNRS), 
        Universit\'e Pierre et Marie Curie (Paris 6), 4 Place Jussieu, 
        75252 Paris, France. 
        \newline
        E-mail: {\tt pgLeFloch@gmail.com}.}  To appear in: {\it Portugal. Math. (2009).}}

\begin{document}

\maketitle


\begin{abstract}
We consider nonclassical entropy solutions to scalar conservation laws with concave-convex
flux functions, whose  set of left- and right-hand admissible states $u_l, u_r$ across undercompressive
shocks is selected by a kinetic function $u_r = \phib(u_l)$.  We introduce 
a new definition for the (generalized) strength of classical and nonclassical shocks, 
allowing us to propose a generalized notion of total variation functional. Relying only upon the natural 
assumption that the composite 
function $\phib \circ \phib $ is uniformly contracting, we prove that the generalized total variation of
front-tracking approximations is non-increasing in time, and we conclude with the existence of nonclassical
solutions to the initial-value problem. We also propose two definitions 
of generalized interaction potentials which are adapted to handle nonclassical entropy solutions
and we investigate their monotonicity properties. In particular, we exhibit an interaction functional 
which is globally non-increasing along a splitting-merging interaction pattern.
\end{abstract}



\pagestyle{myheadings}
\thispagestyle{plain}
\markboth{M. Laforest and P.G. LeFloch}{Diminishing functionals for nonclassical entropy solutions}


\section{Introduction} 
\label{sec:preliminaries}

Consider the following initial-value problem associated with a conservation law in
one-space variable :
\be 
\label{conserv_law}
\aligned
	u_t + f(u)_x & = 0,     
	\\
	u(0,\cdot)  & = u_0,  
\endaligned 
\ee 
where $u_0$ is a function of bounded variation on $\mathbb{R}$, and the
(smooth) flux $f: \RR \to \RR$ is a {\sl concave-convex function,} in the sense that 
\be 
\label{conditions_flux}
\begin{array}{lc} 
 & u \, f''(u) > 0 \quad (u \neq 0), \quad \qquad f'''(0) \neq 0, \\
  &      \displaystyle   \lim_{|u|\to + \infty} f'(u)  = + \infty.
\end{array}          
\ee

Following LeFloch \cite{LeFloch-book}, we consider nonclassical entropy solutions to this problem.
Recall that, in many applications, only a single  entropy inequality can be imposed on the solutions, 
i.e. 
\be 
\label{entropy_ineq}
       U(u)_t + F(u)_x  \leq 0,
\ee 
where the so-called entropy $U$ is a given, strictly convex function and 
the entropy flux $F(u) := \int^u U'(v) f'(v) \, dv$ is determined by $U$. 
It is not difficult to construct multiple weak solutions to the same initial-value problem 
\eqref{conserv_law}, \eqref{entropy_ineq}, so that one realizes that the single
entropy inequality is too lax to determine a unique weak solution. 
In fact, for initial data restricted to lie in one region of concavity or convexity, 
the classical theory applies and leads to a unique entropy solution. 
Non-uniqueness arises when weak solutions contain transitions from positive to negative values, 
or vice-versa.  

The above non-uniqueness property is closely related to the fact that discontinuous solutions,
 in general, depend upon their regularization, that is, 
different regularizations or approximations to the conservation law \eqref{conserv_law}
 may well lead to different solutions
in the limit. This, in particular, is true for solutions to the Riemann problem, corresponding to the special
initial data
\be 
\label{riem} 
u_0(x) = \begin{cases}
u_l, \quad x < 0, 
\\
u_r, \quad x > 0, 
\end{cases}
\ee 
where $u_l, u_r$ are constant states.  Indeed, for a wide class
of regularizations, including regularizations by nonlinear diffusion-dispersion
terms, nonclassical shocks violating Oleinik's entropy inequalities exist. 
The selection of a nonclassical solution is (essentially) equivalent to 
prescribing a {\sl kinetic function} $\phib: \RR \to \RR$  
which, by definition, provides a characterization of admissible nonclassical shocks
connecting states $u_-, u_+$, that is,  
\be 
\label{defn_phib}
    u_+ = \phib(u_-).
\ee

For scalar conservation laws and, more generally, nonlinear hyperbolic systems, 
LeFloch and co-authors initiated the development of a theory of nonclassical entropy solutions 
selected by the kinetic relation \eqref{defn_phib}; various analytical and numerical aspects have been covered. 
We refer to \cite{LeFloch-book} for a review of the theory
and to \cite{LeFlochMohammadian} for recent results on the numerical issues. 
On the other hand, the kinetic relation was originally introduced  
in the context of a hyperbolic-elliptic model describing 
the dynamics of phase transitions in liquids or solids, for which we refer the reader 
to
Slemrod \cite{Slemrod}, 
Truskinovsky \cite{Truskinovsky}, and Abeyaratne and Knowles \cite{AbeyaratneKnowles}.  
LeFloch \cite{LeFloch} used the Glimm scheme and first investigated mathematically the kinetic relations arising in phase dynamics. 
The subject has developed extensively since then, and we will not try here 
to review the literature.

Our objective in the present paper is to discuss the construction of new functionals 
for the (generalized) total variation and wave interaction potential of nonclassical entropy 
solutions to \eqref{conserv_law}. 
We consider solutions generated by Dafermos' front-tracking method, when the 
local Riemann solutions are nonclassical and are determined by a given kinetic relation \eqref{defn_phib}. 
We are interested in deriving uniform estimates for the total variation of solutions and showing that 
the scheme converges to global-in-time, nonclassical entropy solutions 
to the initial-value problem \eqref{conserv_law}.

The new definition of generalized total variation proposed here 
is completely natural, and uses in a direct manner a key property 
satisfied by the kinetic function, 
that is, the uniform contraction property of the second iterate $\phib \circ \phib$. 
Indeed, this property  (see \eqref{contract}, below) motivates our definitions of 
shock strength and, hence, of total variation.
In turn, the proposed construction provides both an improvement of earlier results 
on the subject and a drastic simplification of earlier arguments 
(cf., for instance, with Section~8.1 in \cite{LeFloch-book}). 
Most importantly, 
it appears that our new arguments are robust and may be generalized to tackle systems of equations.

An outline of this paper follows. 
In Section~\ref{sec:prelim}, we begin with a brief review of the theory of kinetic 
relations, by emphasizing what we will need in the rest of this paper. 
We then 
introduce our new definition of generalized shock strength. 
In Section~\ref{sec:TV}, we prove that the proposed generalized total variation functional is 
non-increasing along 
a sequence of Dafermos' front tracking solutions; cf.~Theorem \ref{thm:main}. 
In Section \ref{sec:quad}, we turn to the construction of interaction functionals
adapted to nonclassical entropy solutions. In Theorem~\ref{thm:quad} 
we show that a rather natural definition 
leads to an interaction functional which is non-increasing in all but four interaction patterns. 
Next, in Section \ref{sec:splitmerge} we prove that the proposed interaction functional 
is actually globally non-increasing,
at least in the significant case of merging-splitting wave patterns. 
Such solutions were originally introduced by LeFloch and Shearer \cite{LeFlochShearer}. 
Hence, Theorem~\ref{thm:local_conserv} below demonstrates the 
relevance of the proposed functional to handle nonclassical entropy solutions. 
Finally, in Section~\ref{sec:quad2} we investigate a second definition of interaction functional
which may be better suited for the extension to systems; see Theorem~\ref{thm:quad}. 

\vskip.5cm 


\section{Kinetic function and generalized wave strength}
\label{sec:prelim}

\subsection{Assumptions on the kinetic function}

The general theory of kinetic functions $\phib$ imposes the following conditions on the kinetic 
function: 

\vskip.15cm 

\begin{itemize}
\item[(A1)] The map $\phib : \RR \to \RR$ is Lipschitz continuous and one-to-one. 

\vskip.15cm 

\item[(A2)]  $\partial_u {\phib} (u) < 0 $  for all  $ u $, and $ \phib(0) = 0$.  

\vskip.15cm 

\item[(A3)] The second iterate $ \phib \circ \phib$ is a strict contraction, i.e. for some $K \in (0,1)$
\be 
\label{contract}
       \big| \phib \circ \phib(u) \big| \leq K |u|, \qquad   u \neq 0.
\ee 
\end{itemize}
These are the sole conditions required on the kinetic functions in the present paper. 
We recall that these assumption are satisfied by 
quite general regularizations of the conservation law \eqref{conserv_law}
and, moreover, natural analogues of the properties (A1), (A2), and (A3) 
are known to hold for {\sl systems,} when $\phib$ is properly defined in this more general setting. 
We refer the reader to the monograph \cite{LeFloch-book} and the references cited therein 
for further information.  

The kinetic function determines, via \eqref {defn_phib},
 the set of admissible states $u_-, u_+$ arising on the left- and the right-hand sides of a discontinuity. 
When $u_- > 0$ ($u_-<0$, respectively),
 nonclassical shocks are required when the amplitude of the shock is large enough and 
 the threshold 
$u_+ < \phis(u_-)$ ($\phis(u_-) < u_+$, resp.) is reached. 
This threshold function $\phis: \RR \to \RR$ 
is defined  as the unique value $\phis(u_-) \notin \big\{u_-, \phib(u_-)\big\}$ such that  
$$ 
      \frac{f(u_-) - f\big( \phib(u_-) \big)}{ u_- - \phib(u_-)} 
      =    \frac{f(u_-) - f\big( \phis(u_-) \big)}{ u_- - \phis(u_-)}.
$$
Geometrically, $\big(\phis(u_-),f\big(\phis(u_-) \big)\big)$ is
the point at which the line between $\big(u_-, f(u_-)\big)$
and $\big(\phib(u_-),f\big(\phib(u_-) \big)\big)$ crosses the graph of $f$ 
and therefore, if for instance $u_->0$, the following inequalities must hold  
\be  
\label{phi_order}
   \phib(u_-) < \phis(u_-) < u_-. 
\ee

Under the concave-convex condition \eqref{conditions_flux}, the nonclassical entropy
solution to the Riemann problem with data \eqref{riem}, $u_l >0$, is given by 
\begin{itemize}
\item[i)]  a shock if $\phis(u_l) \leq u_r$, or 
\item[ii)] a nonclassical shock connecting $u_l$ to $\phib(u_l) $
             followed by a classical shock connecting $\phib(u_l)$ to $u_r$
             if $\phib(u_l) < u_r < \phis(u_l) $, or 
\item[iii)]  a nonclassical shock connecting $u_l$ to $\phib(u_l) $
             followed by a classical rarefaction connecting $\phib(u_l)$ to $u_r$
              if $ u_r \leq \phib(u_l) $.           
\end{itemize}
When $u_l <0$ and $u_r>0$, the nonclassical Riemann solver is similar; see \cite{LeFloch} for further details.

Some remarks on our assumptions (A1)--(A3) above are in order. 
In fact, these properties can be motivated from similar properties satisfied 
by the so-called {\sl zero dissipation kinetic function} $\phibo$, 
which is 
determined from the flux $f$ and the given entropy $U$ 
and is {\sl formally} associated to a dispersion-only regularization of \eqref{conserv_law}.
Interestingly, this special kinetic function satisfies a stronger version of 
\eqref{contract}, namely 
\be 
\label{phibo_identity}
            \phibo \circ \phibo(u) = u.
\ee 
On the other hand, for (non-trivial) regularizations including diffusion terms, the kinetic function 
satisfies the strict inequality 
\be
  \label{less_phibo}
        |\phib(u)| < | \phibo(u) |,  \qquad u \neq 0.
\ee 
In the general theory, \eqref{contract} arises as a consequence of the more fundamental 
identities \eqref{phibo_identity} and \eqref{less_phibo}. It is, therefore, natural that our functional 
involve the function $\phibo$ rather than directly $\phib$, and its definition 
makes direct use of the properties  \eqref{phibo_identity} and \eqref{less_phibo}.

Finally, we observe that for any  bounded set of values $|u| < M$, properties (A2) and (A3) 
imply that there exist lower and upper
Lipschitz constants $\underline{\Lip} (u-\phibSqr)>0$ 
and $\overline{\Lip}(u - \phibSqr)< 1$ such that for $u \neq 0$
\be  
\label{strict_contract}
   \underline{\Lip} (u-\phibSqr)  \, \big| u \big| <   \big| u - \phibSqr (u) \big| < \overline{\Lip}(u - \phibSqr) |u|.
\ee 
In the context of the general theory, we note that only the lower bound is non-trivial
and that it does hold for large classes of diffusion-dispersion relations. On the other hand,
for highly degenerate regularizations, such as the one studied in \cite{LeFlochBedjaoui}, 
this condition might be violated.


\subsection{New notion of generalized wave strength}

We introduce here a new functional, which is equivalent to the total variation and 
can be expressed entirely in terms of the kinetic function $\phib$. 
In order to hope for an extension to systems, 
the formulation of the total variation should
let the property \eqref{contract} appear as naturally as possible.
Intuitively, the following definition of strength is based on the idea of
 mapping all the states to the convex region of the flux ($u > 0$), using $\phibo$.
 
\vskip.15cm 
 
\begin{definition}[Notion of generalized wave strength]
 \label{defn:strength}
To each classical or nonclassical wave $(u_-,u_+)$ one associates the 
generalized strength $\sigma(u_-,u_+)$ defined as follows :
when $u_- \geq 0$, 
\be
\mysigma(u_-, u_+) =   \begin{cases}
                                      | u_- - u_+ |,        & u_+  \geq 0, \\
                                      | u_- - \phibo ( u_+) |,   & u_+ < 0,
                                   \end{cases}
\ee
while for $u_- < 0$, 
\be
\mysigma(u_-, u_+) =   \begin{cases}
                                      | \phibo( u_- ) - \phibo( u_+ ) |,        & u_+  < 0,
                                       \\
                                      | \phibo( u_- ) -  u_+ |,   & u_+ \geq 0.
                                   \end{cases}
\ee
\end{definition}

\vskip.15cm 

The following notation will be helpful for the analysis.
Henceforth, classical shocks and rarefactions joining two positive states is denoted by 
$\Cp$ and $\Rp$, respectively. Similarly, when both neighboring states are negative, we 
 write $\Cm$ and $\Rm$. When a shock joins a positive state with a
negative state, we write $\Cpm$ or $\Npm$ depending on whether or
not the shock is classical or nonclassical, respectively. 
When the signs of the neighboring states
are reversed, we simply write $\Cmp$ and $\Nmp$. 
Shocks $\Cpm$ and $\Cmp$ are also sometimes called crossing shocks.

The generalized wave strength enjoys several immediate properties. 

\begin{itemize}

\item The first important property to observe is that the proposed generalized strength is 
{\sl continuous} as $u_+$ crosses $\phis(u_-)$ and the solution of the Riemann
problem goes from a single crossing shock to a nonclassical shock
followed by a classical shock.  

For $u_- >0$, this is checked from the
inequalities $\phib(u_-) < \phis(u_-) < 0 $ and the property \eqref{contract} which 
yield 
\begin{align*}
     \mysigma(u_-,\phis(u_-)) & = \big| u_- - \phibo \circ \phis( u_-) \big| \\
                                            & = \big| u_l - \phibSqr( u_-) \big| 
                                               +  \big| \phibSqr( u_-) -  \phibo \circ \phis( u_-) \big|  \\
                                             & = \mysigma(u_-, \phib(u_-)) 
                                                 + \mysigma( \phib(u_-), \phis(u_-)).    
\end{align*}
In fact, $\phis(u_-)$ could be positive and a similar computation
would show that continuity still holds.

\item The second important property of the generalized strength is its ``equivalence'' with
the usual notion of strength. For positive non-crossing shocks  and rarefactions,
the definition is the same as the usual one. When the rarefaction
and the non-crossing shocks have two negative neighboring states, then
there exists a positive constant $C'$ such that
\be  
\label{positivity1}
    \mysigma (u_-, u_+) = \big| \phibo (u_-) - \phibo(u_+)  \big| 
                                  \geq C'  | u_- - u_+ |,
\ee 
as long as $u_-, u_+$ stay within a bounded neighbourhood of the origin. 
For the shocks $\Cpm$ and
$\Npm$, it suffices to use  \eqref{strict_contract} to show that
the definition is equivalent to the usual notion of strength :  
\be 
 \label{positivity2}
         \sigma( \Npm ) = \big| u_- - \phibSqr( u_- ) \big| 
                              > C''  2 \big| u_- \big| > C''' \big|  u_-  - \phib(u_-)  \big|. 
\ee 
The same argument also applies to shocks $\Cmp$ and $\Nmp$.
Interestingly, our definition of shock strength makes clear the intuitive idea
that it should be increasingly difficult to measure the strength of nonclassical shocks 
as the zero-diffusion (dispersion-only) case is approached.

\end{itemize}


\vskip.5cm 

\section{Diminishing total variation functional}
\label{sec:TV}

 We now introduce 
 front-tracking approximate solutions to \eqref{conserv_law} based on a 
 nonclassical Riemann solver, as was constructed by Dafermos \cite{Dafermos72}
 in the classical situation. 
 These approximations are piecewise constant 
 in space and are determined from the nonclassical Riemann solver
 described in the previous section.
  
 The first step of their construction is to build a 
 piecewise constant approximation  of 
 the initial data $u_0$ which admits finitely many 
 discontinuities and approaches $u_0$ in the $L^1$ norm with an error $\eps$,
 for some small $\eps$. The Rankine-Hugoniot condition can be 
 used to propagate, in a conservative manner, the discontinuities of the initial data.
 When the Riemann solver calls for continuous waves, we 
 replace them by a sequence of small discontinuities $(u_-, u_+)$  whose strength satisfy $\sigma(u_-,u_+) < \eps$.
 
 When two discontinuities meet, the nonclassical Riemann solver is used to
 determine the neighboring states, but we continue to enforce
 that all outgoing waves be discontinuities. Despite the fact that the number of outgoing waves may be
 larger than the number of incoming discontinuities, one can check 
 (see \cite{LeFloch-book} for details)
 that the total
 number of discontinuities remains bounded for all times, so that
  the front-tracking approximation
 can be defined for all times.

For a front-tracking approximation $u: \RR_+ \times \RR \to \RR $, 
formed entirely of propagating discontinuities (each denoted by $\alpha$),  
the inequalities \eqref{positivity1} and \eqref{positivity2}, as well as the fact that
the kinetic functions are Lipschitz continuous, imply that 
\be 
  \label{myTV}
   V \big(  u(t) \big) := \sum_{\alpha}  \mysigma(  u_-^\alpha, u_+^\alpha ),
\ee  
is equivalent to the total variation norm 
\be 
  \label{TV}
   \operatorname{TV} \big(  u(t) \big) := 
           \sum_\alpha  \big|  u_-^\alpha - u_+^\alpha \big|,
\ee 
where $u_\pm^\alpha$ denote the left- and right-hand states of the discontinuity $\alpha$. 

Incidently,
we mention that  the kinetic function $\phibo$ induces an isometry on the space
of $\operatorname{BV}(\RR)$ functions with respect to the norm $V$ :
\be 
   \label{defn:Phi_operator}
     \begin{array}{rrcl}
      \Phi^{\flat}_0 : & \operatorname{BV}(\RR) & \to & \operatorname{BV}(\RR) \\
                             & u(\cdot)                          & \mapsto & \phibo \circ u(\cdot).
     \end{array}                        
\ee 
Using \eqref{phibo_identity}, it is straightforward to check that $\Phi_0^{\flat}$ is an isometry. 
This result will not be used directly in this paper.


\vskip.15cm 

\begin{theorem}[Diminishing total variation functional]
 \label{thm:main}
Let $\phib$ be a kinetic function satisfying the properties (A1)--(A3). 
For every front-tracking approximation $u: \RR_+ \times \RR \to \RR$ to the conservation law
\eqref{conserv_law} based on the nonclassical Riemann solver associated with $\phib$, 
the generalized total variation functional
$V\big( u(t) \big)$ is non-increasing. 
Precisely, the change in $V$ 
during an interaction is given by: 
\be 
  \big[ V \big] \leq \begin{cases}
                             - c_1 \mysigma( R^{in} ),  & \text{Cases RC-1, 
                               RC-3, CR-1, CR-2, CR-4,} \\
                             - c_2 \mysigma( R^{in} ), & \text{Cases RC-2, RN, } \\
                             - 2 \big(  \mysigma( R^{in} ) - \mysigma( R^{out} )  \big),  & 
                                    \text{Case CR-3, } \\
                              0,  &   \text{all other cases.}
                       \end{cases}
\ee
Here, $c_1 = \min \{ 1, \underline{\Lip}( \phibo) \} $, $c_2 =  \underline{\Lip}( u -\phibSqr )$,
while 
$R^{in}$ and $R^{out}$ denote the incoming and outgoing rarefactions at an interaction, respectively. 
(The list of interactions is specified in the proof below.) 
\end{theorem}

\vskip.15cm 

From this theorem we immediately deduce the following:

\vskip.15cm 

\begin{corollary}[Existence of nonclassical entropy solutions] 
For any initial data $u_0 \in L^{\infty}(\RR) \cap \operatorname{BV}(\RR) $
of \eqref{conserv_law}, and any sequence of front-tracking approximations $u^{\eps}$
such that $u^{\eps}(\cdot,0)$ converges to $u_0(\cdot)$ in $L^1(\RR)$,
there exists a subsequence of front-tracking approximations that converge
in $\Lip \big( [0,T), L^1(\RR) \big) \cap L^{\infty}([0,T),\operatorname{BV}(\RR))$
to a solution of the initial value problem \eqref{conserv_law}. 
\end{corollary}

\vskip.15cm 

{\it Proof of Theorem \ref{thm:main}.}  
Our proof is a generalization of the one given in \cite{LeFloch-book}
(see Section~8.1).
We need here to compute the variation of our functional $V$ by distinguishing between
$16$ possible interactions, after assuming $u_l >0$ for definiteness, 
as described in Section~4.3 of \cite{LeFloch-book}. 
We note that we have one more case to consider here, 
 since 
our assumptions on the kinetic function are general enough to  
allow for the so-called CC-3 interactions. 
The subscript $'$ is used to indicate that a wave is outgoing and with some abuse of
 notation for the wave strength we write, for example,
$\mysigma( \Npm)$ for the shock strength of a nonclassical wave.
During a generic interaction between two waves, we denote
the states on both sides of the left-hand wave by $u_l$ and $u_m$ while
those associated with the right-hand wave are denoted by $u_m$ and $u_r$. 
Finally, the bounds on the change of $V$ depend only on the
properties \eqref{contract}--\eqref{strict_contract} which are therefore used freely throughout.

\vskip.15cm 

\begin{flushleft}
\textbf{Case RC-1 :} ($\Rp \CDown$)-(${\CDown}'$). 
This case is determined by the constraints
$$
    \max \big( \phis(u_l), \phis(u_m) \big) < u_r < u_l, \qquad
      0 < u_l < u_m.
$$
This is further subdivided into two subcases depending on the sign of $u_r$.
When $u_r>0$, then the interactions are entirely
classical ($\Rp \Cp$)-(${\Cp}'$) and the inequalities
$ 0 < u_r < u_l < u_m $ suffice to check that
\begin{equation*}
  [V] = \mysigma\big( {\Cp}'\big) 
        - \mysigma\big( \Rp \big) -\mysigma\big( \Cp \big)
        = - 2 |u_m - u_l| = - 2 \mysigma\big( \Rp \big).
\end{equation*}
When $u_r<0$, then the interaction involves crossing shocks
($\Rp \Cpm$)-(${\Cpm}'$) and the states involved in measuring the 
strengths of the waves are
\be  \label{RC1_SUB2_states}
        0 < \phibo(u_r) < u_l < u_m,
\ee
since $\phibo(u_r) < \phibSqr(u_l) < u_l$. Then
\begin{align*}
  [V] & = \mysigma\big( {\Cpm}'\big) 
        - \mysigma\big( \Rp \big) -\mysigma\big( \Cpm \big) \\
       & = |u_l - \phibo(u_r)| - |u_m - u_l| - |u_m - \phibo(u_r) | 
        = - 2 \mysigma\big( \Rp \big).
\end{align*} 

\vskip.15cm 

\textbf{Case RC-2 :} ($\Rp \Cpm$)-(${\Npm}' {\Rm}'$). 
This case is defined by
$$
  \phis(u_m) < u_r \leq \phib(u_l) < 0 < u_l < u_m.
$$
In the first subcase, we assume $\phibo(u_r) < u_l$ and use the
previous conditions to deduce 
\be  \label{RC2_SUB1_states}
        \phibSqr(u_l) < \phibo(u_r) < u_l < u_m.
\ee
The analysis when $u_l < \phibo(u_r)$ will not be treated
since only the order of those two terms changes 
but the conclusions remain the same.
\begin{align*}
  [V] & =    \mysigma\big( {\Npm}' \big) + \mysigma\big( {\Rm}'\big) 
              - \mysigma\big( \Rp \big) -\mysigma\big( \Cpm \big) \\
       & =   |u_l - \phibSqr(u_l)| + | \phibo(u_r) - \phibSqr(u_l) |
               - |u_m-u_l | - |u_m - \phibo(u_r)|  \\
       & = 2 | \phibo(u_r) - \phibSqr(u_l) | - 2 | u_m - u_l |.        
\end{align*}
Since $\phis(u_m) < u_r <0$, from property \eqref{contract} we deduce
$0 < \phibo(u_r) < \phibo \circ \phis(u_m) < \phibSqr(u_m)$. 
Combining this with \eqref{RC2_SUB1_states}, we find
$$
 [V] \leq - 2 \underline{\Lip}( u - \phibSqr ) |u_m - u_l| = -2 \underline{\Lip}( u - \phibSqr ) \mysigma( \Rp ).
$$

\vskip.15cm 

\textbf{Case RC-3 :} ($\Rp \CDown$)-(${\Npm}' {\CUp}'$). 
The conditions initially satisfied by the neighboring states of the incoming waves  are
$$
   \max\big( \phib(u_l), \phis(u_m) \big) < u_r < \phis(u_l), \qquad 0 < u_l < u_m. 
$$
In the first subcase, we assume $u_r <0$ and the
interaction is ($\Rp \Cpm$)-(${\Npm}' {\Cm}'$)
with the following states appearing in the strength of the waves
\be  \label{RC3_SUB1_states}
         \phibo(u_r) < \phibSqr(u_l) < u_l < u_m.
\ee
Therefore, the change in $V$ is
\begin{align*}
  [V] & =    \mysigma\big( {\Npm}' \big) + \mysigma\big( {\Cm}'\big) 
              - \mysigma\big( \Rp \big) -\mysigma\big( \Cpm \big) \\
       & = | u_l - \phibSqr(u_l) | + |\phibSqr(u_l) - \phibo(u_r)| 
             - |u_m-u_l| - | u_m - \phibo(u_r) | \\
       & = - 2 | u_m - u_l |  = - 2 \mysigma( \Rp ).        
\end{align*}
In the second subcase with $0< u_r$, we have
$u_r < \phis(u_l) < \phib \circ \phib (u_l) < \phibo \circ \phib (u_l)$
and therefore
\be  \label{RC3_SUB2_states}
         u_r < \phibSqr(u_l) < u_l < u_m.
\ee
The change in $V$ is then
\begin{align*}
  [V] & =    \mysigma\big( {\Npm}' \big) + \mysigma\big( {\Cmp}'\big) 
              - \mysigma\big( \Rp \big) -\mysigma\big( \Cp \big) \\
       & = | u_l - \phibSqr(u_l) | + |\phibSqr(u_l) - u_r| 
             - |u_m-u_l| - | u_m - u_r | \\
       & = - 2 | u_m - u_l |  = - 2 \mysigma( \Rp ).        
\end{align*}

\vskip.15cm 

\textbf{Case RN :} ($\Rp \Npm$)-(${\Npm}' {\Rm}'$). 
The states on both sides of the waves satisfy
$$
  0 < u_l < u_m \qquad \text{and} \qquad u_r = \phib(u_m).
$$
Two cases again occur depending on the relative order of $u_l$
with respect to $\phibo(u_r) = \phibSqr(u_m)$. We consider only the case where
\be  \label{RN_SUB1_states}
         \phibSqr(u_l) < \phibSqr(u_m) < u_l < u_m,
\ee
the other, $\phibSqr(u_l) < u_l < \phibSqr(u_m)  < u_m $, being similar.
In this case, the change in $V$ is
\begin{align*}
  [V] & =    \mysigma\big( {\Npm}' \big) + \mysigma\big( {\Rm}'\big) 
              - \mysigma\big( \Rp \big) -\mysigma\big( \Npm \big) \\
       & = | u_l - \phibSqr(u_l) | + | \phibSqr(u_m) - \phibSqr(u_l) | \\
       & \phantom{ = }      - |u_m-u_l| - | u_m - \phibSqr(u_m) | \\
       & =  2 |\phibSqr(u_m) - \phibSqr(u_l) | - 2| u_m - u_l |  \\
       & \leq -2 \underline{\Lip}( u - \phibSqr) | u_m - u_l | 
           = -2 \underline{\Lip}( u - \phibSqr) \mysigma( \Rp ).        
\end{align*}

\vskip.15cm 

\textbf{Case CR-1 :} ($\Cpm \Rm$)-(${\Cpm}'$).
This is a simple case where the states are initially ordered as
$$
        \phis(u_l) < u_r < u_m \leq 0 < u_l.
$$
This provides $\phibo(u_m) < \phibo(u_r) < \phibo( \phis(u_l) ) < u_l$. The result is :
\begin{align*}
  [V] & =   \mysigma\big( {\Cpm}'\big) 
          - \mysigma\big( \Cpm \big) -\mysigma\big( \Rm \big) \\
       & = |u_l - \phibo(u_r) | - |u_l - \phibo(u_m) | - |\phibo(u_r) - \phibo(u_m)| \\
       & = - 2 |\phibo(u_r) - \phibo(u_m)| \leq -2 \underline{\Lip}( \phibo ) \mysigma( \Rm ). 
\end{align*}

\vskip.15cm 

\textbf{Case CR-2 :} ($\Cp \Rp$)-(${\Cp}'$). 
The waves are entirely classical since
$0 \leq u_m < u_r < u_l$. There is nothing new to check
and the change in $V$ is immediately found to be 
$$
 [V] = - 2 |u_r - u_m| = -2 \mysigma( \Rp ).
$$

\vskip.15cm 

\textbf{Case CR-3 :} ($\Cpm \Rm$)-(${\Npm}' {\Rm}'$).
The states begin in the order
$$
  u_r \leq \phib(u_l) < \phis(u_l) < u_m \leq 0 < u_l.
$$
One subcase is obtained when we assume that $u_l < \phibo(u_r)$.
The important states appearing in the definition of the
strength of waves are then ordered as
\be
\label{CR3_SUB1_states}
         \phibo(u_m) < \phibSqr(u_l) < u_l < \phibo(u_r).
\ee
The change in $V$ is then found to be
\begin{align*}
  [V] & =    \mysigma\big( {\Npm}' \big) + \mysigma\big({\Rm}' \big) 
              - \mysigma\big( \Cpm \big) -\mysigma\big( \Rm \big) \\
       & = | u_l - \phibSqr(u_l) | + | \phibo(u_r) - \phibSqr(u_l) | 
             - | u_l - \phibo(u_m) | - |\phibo(u_r) - \phibo(u_m) | \\
        & =  - 2 | \phibSqr(u_l) - \phibo(u_m) |.      
\end{align*}
In the second subcase, when $\phibo(u_r) < u_l$, we use
\be
\label{CR3_SUB2_states}
         \phibo(u_m) < \phibSqr(u_l) < \phibo(u_r) < u_l 
\ee
to deduce the same equality :
\begin{align*}
  [V] & = | u_l - \phibSqr(u_l) | + | \phibo(u_r) - \phibSqr(u_l) | 
             - | u_l - \phibo(u_m) | - |\phibo(u_r) - \phibo(u_m) | \\
        & =  - 2 | \phibSqr(u_l) - \phibo(u_m) |.      
\end{align*}
We now observe that in both subcases, the change is equal to 
a physically relevant quantitiy  
$$
  [V]  = -2 \big(  \mysigma( \Rm ) - \mysigma( {\Rm}' ) \big).
$$

\vskip.15cm 

\textbf{Case CR-4 :} ($\Cpm \Rm$)-(${\Npm}' {\Cm}'$).
The states of the incoming waves satisfy
\be 
     \phib(u_l) < u_r < \phis(u_l) < u_m \leq 0 < u_l.
\ee
This leads us to the inequalities
\be
\label{CR4_states}
        \phibo(u_m) < \phibo(u_r)  < \phibSqr(u_l) < u_l.
\ee
The change in our functional is therefore
\begin{align*}
  [V] & =    \mysigma\big( {\Npm}' \big) + \mysigma\big({\Cm}' \big) 
              - \mysigma\big( \Cpm \big) -\mysigma\big( \Rm \big) \\
       & = | u_l - \phibSqr(u_l) | + | \phibSqr(u_l) - \phibo(u_r) | 
             - | u_l - \phibo(u_m) | - |\phibo(u_r) - \phibo(u_m) | \\
        & =  - 2 | \phibo(u_r) - \phibo(u_m) | \leq -2 \underline{\Lip}(\phibo) \mysigma( \Rm ).      
\end{align*}

\vskip.15cm 

\textbf{Case CC-1 :} ($\Cp \CDown$)-(${\CDown}'$).
This is another simple case. We begin with
$$
     \max\big( \phis(u_l), \phis(u_m)  \big) < u_r < u_m < u_l,
     \qquad \text{and} \qquad 0 \leq u_m. 
$$
When $0 \leq u_r$, all the waves are classical and it is easy to show
that $[V] = 0$. When $u_r < 0$, we still have
$ \phib(u_m) < u_r$ and therefore the important states are ordered
$$
         \phibo(u_r) < u_m < u_l.
$$
It then easy to check that even in this case, $[V] =0$.

\vskip.15cm 

\textbf{Case CC-2 :} ($\Cpm \CUp$)-(${\CDown}'$).
This interaction is constrained by the initial states in the following manner
$$
    \phis(u_l) < u_m < u_r < \phis(u_m) < u_l, 
      \qquad \text{and} \qquad u_m < 0. 
$$
Two subcases appear depending on the sign of $u_r$. When $u_r >0$,
the interaction is ($\Cpm \Cmp$)-(${\Cp}'$). Since $\phib(u_l) < u_m$,
$\phibo(u_m) < \phibSqr(u_l) < u_l$. Combining this with
the fact that $u_r < \phis(u_m) < \phibo(u_m)$, we can deduce
that the states appearing in $[V]$ are ordered
$$
        u_r < \phibo(u_m) < u_l. 
$$
The change in $V$ is then
\begin{align*}
  [V] & =    \mysigma\big( {\Cp}' \big)  
              - \mysigma\big( \Cpm \big) -\mysigma\big( \Cmp \big) \\
       & = | u_l - u_r | 
             - | u_l - \phibo(u_m) | - |\phibo(u_m) - u_r | = 0.
\end{align*}
The second subcase treats $u_r <0$ and interactions
($\Cpm \Cm$)-(${\Cpm}'$). The states used in our definition of the 
strength of the waves are 
$$
       \phibo(u_r) < \phibo(u_m) < u_l. 
$$
A short computation shows that
\begin{align*}
  [V] & =    \mysigma\big( {\Cpm}' \big)  
              - \mysigma\big( \Cpm \big) -\mysigma\big( \Cm \big) \\
       & = | u_l - \phibo(u_r) | 
             - | u_l - \phibo(u_m) | - |\phibo(u_m) - \phibo(u_r) | = 0. \\
\end{align*}

\vskip.15cm 
 
\textbf{Case CC-3 :} ($\Cp \CDown$)-(${\Npm}' {\CUp}'$ ).
This interaction represents a typical transition from one
crossing shock to a nonclassical shock. The states are
$$
  \phib(u_l) < \phis(u_m) < u_r < \phis(u_l) < u_m < u_l,
\qquad 0 \leq u_m. 
$$
The first (and most common) subcase occurs when $u_r < 0$
and the interaction is ($\Cp \Cpm$)-(${\Npm}' {\Cm}'$ ).
The first subcase needs to further subdivided into
two cases. When $u_m < \phibSqr(u_l)$, all the important states
are ordered
\be 
\label{CC3_SUB1.1_states}
       \phibo(u_r) < u_m < \phibSqr(u_l) < u_l.
\ee
Under these circumstances, the functional $V$ doesn't change. 
\begin{align*}
  [V] & =    \mysigma\big( {\Npm}' \big) + \mysigma\big({\Cm}' \big) 
              - \mysigma\big( \Cp \big) -\mysigma\big( \Cpm \big) \\
       & = | u_l - \phibSqr(u_l) | + | \phibSqr(u_l) - \phibo(u_r) | 
             - | u_l - u_m | - | u_m - \phibo(u_r) | = 0.
\end{align*}
On the other hand, when $\phibSqr(u_l) < u_m$ in the first subcase,
we have
\be
\label{CC3_SUB1.2_states}
       \phibo(u_r) <  \phibSqr(u_l) < u_m < u_l.
\ee
It is easy to see that we again have $[V] = 0$. 

In the second subcase, $0 < u_r$, we again need to introduce
two additional subcases to handle the interactions
($\Cp \Cp$)-(${\Npm}' {\Cmp}'$ ). When $u_m < \phibSqr(u_l)$, then the
states used in the definition of wave strengths are
$$
       u_r < u_m < \phibSqr(u_l) < u_l.
$$
The change in $[V]$ is
\begin{align*}
  [V] & =    \mysigma\big( {\Npm}' \big) + \mysigma\big({\Cmp}' \big) 
              - \mysigma\big( \Cp \big) -\mysigma\big( \Cp \big) \\
       & = | u_l - \phibSqr(u_l) | + | \phibSqr(u_l) - u_r | 
             - | u_l - u_m | - | u_m - u_r | = 0.
\end{align*}
When $0<u_r$ and $\phibSqr(u_l) < u_m$, we obtain the same result. 

\vskip.15cm 

\textbf{Case CN-1 :} ($\Cp \Npm$)-(${\Cpm}'$).
The states defining the waves are characterized by the inequalities
$$
   0 < u_m < u_l \qquad \text{and} \qquad 
      \phis(u_l) \leq u_r = \phib(u_m).
$$
This implies that
$$
       \phibSqr( u_m ) < u_m  < u_l.
$$
We deduce
\begin{align*}
  [V] & =    \mysigma\big({\Cpm}' \big) 
              - \mysigma\big( \Cp \big) -\mysigma\big( \Npm \big) \\
       & = | u_l - \phibSqr(u_m) | 
             - | u_l - u_m | - | u_m - \phibSqr (u_m) | = 0.
\end{align*}

\vskip.15cm 

\textbf{Case CN-2 :} ($\Cpm \Nmp$)-(${\Cp}'$).
We begin with states satisfying
$$
    \phis(u_l) < u_m < 0 \qquad \text{and} \qquad
       u_r = \phib(u_m).
$$
To measure the wave strengths, we observe that
$\phib(u_l) < u_m$ implies $\phibo(u_m) < \phibSqr(u_l) < u_l$
and therefore 
$$
        u_r = \phib(u_m) < \phibo(u_m) < u_l.
$$
The jump in $V$ is now
\begin{align*}
  [V] & =    \mysigma\big({\Cp}' \big) 
              - \mysigma\big( \Cpm \big) -\mysigma\big( \Nmp \big) \\
       & = | u_l - \phib(u_m) | 
             - | u_l - \phibo(u_m) | - | \phibo(u_m) - \phib(u_m) | = 0.
\end{align*}

\vskip.15cm 

\textbf{Case CN-3 :} ($\Cp \Npm$)-(${\Npm}' {\Cm}'$).
The states are initially ordered as
$$
      0 < u_m < u_l \qquad \text{and} \qquad 
       u_r = \phib(u_m) < \phis(u_l).
$$
A first subcase occurs when $u_m < \phibSqr(u_l)$.
To compute the change in $V$, we can then use the inequalities
\be \label{CN3_SUB1_states}
       \phibSqr( u_m ) < u_m  < \phibSqr(u_l) < u_l,
\ee
to deduce
\begin{align*}
  [V] & =    \mysigma\big( {\Npm}' \big) + \mysigma\big({\Cm}' \big) 
              - \mysigma\big( \Cp \big) -\mysigma\big( \Npm \big) \\
       & = | u_l - \phibSqr(u_l) | + | \phibSqr(u_l) - \phibSqr(u_m) |  \\
       & \phantom{ =}   - | u_l - u_m | - | u_m - \phibSqr(u_m) | \\
       &   = 0.
\end{align*}

Similarly, if $\phibSqr(u_l) < u_m$, then the important inequalities become 
\be
    \label{CN3_SUB2_states}
       \phibSqr( u_m ) < \phibSqr(u_l) < u_m < u_l,
\ee
and $[V] = 0$.

\vskip.15cm 

\textbf{Case NC :} ($\Npm \CUp$)-(${\CDown}'$).
This interaction is constrained by the states
$$
  u_m = \phib(u_l) \qquad \text{and} \qquad
       \phis(u_l) < u_r < \phis(u_m) < u_l.
$$
The first subcase occurs when $u_r < 0$ and the interaction
is precisely ($\Npm \Cm$)-(${\Cpm}'$). The important states are then
ordered as
$$   
       \phibo( u_r ) < \phibSqr(u_l) < u_l.
$$
With these observations, we find that
\begin{align*}
  [V] & =    \mysigma\big( {\Cpm}' \big) 
              - \mysigma\big( \Npm \big) -\mysigma\big( \Cm \big) \\
       & = | u_l - \phibo(u_r) | 
             - | u_l - \phibSqr(u_l) | - |\phibSqr(u_l)  - \phibo(u_r) | = 0.
\end{align*}
In the second subcase, $u_r > 0$, it is easy to check that
$$ 
       u_r  < \phibSqr(u_l) < u_l,
$$
and $[V] = 0$.

\vskip.15cm 

\textbf{Case NN :} ($\Npm \Nmp$)-(${\Cp}'$).
This interaction is the limiting case $u_r \to \phib(u_m)$
of Case NC. By continuity of wave strengths, we must
also have $[V]=0$. 

\end{flushleft}
\hfill $\Box$


\section{Quadratic interaction potential (part 1)}
\label{sec:quad}

In the rest of this paper, we investigate two quadratic interaction potentials, keeping in mind 
from experience with classical shock waves,
that different functionals may be of particular interest in different circumstances. 
We begin by searching for a functional of the form
\be  
\label{defn:quad}
Q\big( u(t) \big) := \sum_{\text{$\alpha$ approaches $\beta$}}
                   \mysigma( u^{\alpha}_l, u^{\alpha}_r) \mysigma ( u^{\beta}_l, u^{\beta}_r ),  
\ee
where the proper definition of ``pairs of approaching waves" is essential and is now specified.

In Glimm's original paper \cite{Glimm} for systems of conservation laws, a
definition was proposed which, in the scalar case, amounted to stating that
pairs of waves were always approaching unless both are rarefactions.
The purpose of this section is to investigate the original definition of Glimm in the context of 
nonclassical shocks.  

\vskip.15cm 

\begin{definition}  \label{defn:approaching}
A wave $\alpha$ is said to {\sl weakly approach} a wave $\beta$, unless both are rarefaction waves. 
As far as this definition is concerned, waves $\Rp, \Rm$ are both to be considered as rarefaction waves
and $\Cp, \Cpm, \Npm, \Cm, \Cmp, \Nmp$ are all to be considered as shock waves.
\end{definition}

\vskip.15cm 

Our main result in the present section is as follows. 

\vskip.15cm 

\begin{theorem}[``Weak interaction'' potential for nonclassical shocks]
\label{thm:quad}
Let $\phib$ be a kinetic function satisfying the properties (A1)--(A3). 
Consider the functional $Q_{\text{weak}}$ defined by \eqref{defn:quad} where the summation is made over all 
weakly interacting waves in the sense of Definition~\ref{defn:approaching}. Then, when evaluated on a sequence of front-tracking solutions, 
$Q_{\text{weak}}$ is strictly decreasing during all interactions except in the cases 
RC-3, CR-4, CC-3 and CN-3. In fact, for each of these exceptional interactions, there 
exist initial data for which $V + C_0 \,  Q_{\text{weak}}$ is increasing during the interaction 
for every positive $C_0$. 
\end{theorem}

\vskip.15cm 
 
In contrast, in Chapter 8 of \cite{LeFloch-book} in joint work by Baiti, LeFloch, and Piccoli, 
different definitions of both wave strengths and approaching waves are used 
and the resulting Glimm functional $V + K \, Q$ is strictly decreasing for some $K$. 
In this sense, the interaction functional $Q_{\text{weak}}$ above  
may appear to be less satisfactory. However, our assumptions on the kinetic function are {\sl 
completely natural} --a major advantage toward a future extension to systems--  
and, furthermore, an analysis 
of ``splitting/merging'' solutions (in the following section) will show that {\sl globally
in time} the functional $Q_{\text{weak}}$ does decrease. 
 
 \vskip.15cm 

Several justifications for our definition of potential are now provided, 
the strongest argument being the requirement of continuity: 
\begin{enumerate} 

\item[1.] Given that $V$ is continuous in 
$\operatorname{BV}(\RR)$  and that  $\Phi^{\flat}_0$ is already an isometry
with respect to this $\operatorname{BV}(\RR)$ norm, it is tempting to assume that any reasonable 
interaction potential $Q$ should  also be continuous in $\operatorname{BV}(\RR)$.
We observe that any shock $\Cp$ can be 
continuously deformed
(as measured by Definition \ref{defn:strength})
 into, first, a
crossing shock $\Cpm$ and, then, a pair of shocks $\Npm \Cm$. 
It is easy to see that imposing continuity would imply Definition \ref{defn:approaching}.

\item[2.]  Another argument can be found by looking at a class of solutions called
splitting/merging solutions, introduced in \cite{LeFlochShearer} and discussed further 
in Section \ref{sec:splitmerge}. These solutions illustrate that some initial
data can go through a nearly periodic process of creation and destruction of nonclassical
shocks. In particular, nonclassical shocks can indirectly have non-trivial
interactions with shocks on their right-hand side. 

\item[3.] In \cite{LeFloch-book}, nonclassical shocks are precluded from interacting with their right-hand neighbours, 
and 
it is argued that nonclassical shocks are (slow) undercompressive and, thus, move away from their
right-hand neighbors. However, this definition of approaching waves
ignores the above-mentioned possibility of 
nonclassical shocks having indirect interactions with shocks on their right-hand side. 
In any case, such an interaction functional then would not be continuous in $\operatorname{BV}(\mathbb{R})$.

\end{enumerate}

\vskip.15cm 

For the proof of Theorem~\ref{thm:quad} we proceed as follows. 
During any isolated interaction between two waves in an approximation, 
the change in $Q= Q_{\text{weak}}$ is of two types
\be
\label{decompo} 
   [Q] = [Q]_1 + [Q]_2.
\ee
In this decomposition, $[Q]_1$ denotes the change in the products of the strengths of waves either 
incoming or exiting the interaction and $[Q]_2$ denotes the change in products of strengths
of waves where only {\sl one} of the waves was directly involved in the interaction.
Moreover, if a wave $C$ is involved in an interaction, we define
$$
    W(C) :=   \sum_{\substack{\text{$B$ approaches $C$}\\ \text{$B$ did not interact}}}
                    \mysigma ( B ).
$$ 
According to Definition \ref{defn:approaching}, if the incoming wave $C_{\text{in}}$ and
the outgoing wave $C_{\text{out}}$ are of the same type (i.e. both of rarefaction or shock type), then
$W(C_{\text{in}}) = W(C_{\text{out}})$. Theorem \ref{thm:quad} requires the following lemma 
which emphasize the (wrong) sign of the terms $[Q]_1$ appearing in the most difficult interactions.
The lemma gives a precise characterization of the failings of $Q_{\text{weak}}$ for those interactions.

\vskip.15cm 

\begin{lemma}  
\label{lemma:quad}
The interactions involving a crossing shock, classical or not, and a small wave $W_{\text{in}}$ are 
RC-3, RN, CR-4, CC-3 and CN-3. Suppose that the states neighboring the two consecutive waves are
$u_l, u_m$ and $u_r$. For $u_l$ fixed and as a function of the strength of the incoming small wave, 
the change in $[Q]_1$ is as follows :
\begin{itemize}
\item[ ]RN : $[Q]_1 < 0$ for all $R_{\text{in}}$;
\item[] CN-3 : $[Q]_1 < 0$ if  $C_{\text{in}}$ is sufficiently weak;
\item[] RC-3 : $[Q]_1 < 0$ if  $R_{\text{in}}$ is sufficiently strong;
\item[] CR-4, CC-3 : $[Q]_1$ can take on both signs, depending on $u_m$ and $u_r$.
\end{itemize}
Moreover, for RC-3 and CR-4 interactions, the largest positive change occurs when $\sigma(R_{\text{in}})=0$.
\end{lemma}

\

{\it Proof.} We distinguish between several cases. 

\begin{flushleft}

\textbf{Case RN :} ($\Rp \Npm$)-(${\Npm}' {\Rm}'$).
The states appearing in the shock strengths describe two subcases:
\begin{align*}
      \text{ either } \,  & \phibSqr(u_l) < \phibSqr(u_m) < u_l < u_m, 
       \\
        \text{ or } \,  &  \phibSqr(u_l) < u_l < \phibSqr(u_m) < u_m. 
\end{align*}
For both sets of inequalities $  \mysigma({\Rm}') < \mysigma(\Rp)$ and 
$  \mysigma({\Npm}') < \mysigma(\Npm) $. As a result,
$$
     [Q]_1 =    \mysigma({\Npm}') \mysigma({\Rm}')  -  \mysigma(\Rp) \mysigma( \Npm)  < 0.
$$

\textbf{Case RC-3 :} ($\Rp \CDown$)-(${\Npm}' {\CUp}'$).
The fundamental states describe two subcases \eqref{RC3_SUB1_states} and \eqref{RC3_SUB2_states},
here rewritten
\begin{align*}
     \mbox{ either } \,  & \phibo(u_r) < \phibSqr(u_l) < u_l < u_m, 
      \\
    \mbox{ or } \,  &  u_r < \phibSqr(u_l) < u_l < u_m.
\end{align*}
Focusing on the first case, the change in $[Q]_1$ is
\begin{align*}
      [Q]_1(u_l, u_m, u_r) = &  \mysigma({\Npm}') \mysigma( {\CUp}' ) - \mysigma( \Rp )\mysigma(\Cpm) \\
                 = & \big( u_l-\phibSqr(u_l) \big) \, \big( \phibSqr(u_l) - \phibo(u_r) \big) 
                       - (u_m - u_l) \, \big( u_m - \phibo(u_r) \big). 
\end{align*}
Clearly, the RC-3 interaction can still occur even as the strength $\mysigma(\Rp) \to 0$,
therefore the expression above shows that $[Q]_1 > 0$ in that limiting case. In fact, we will show
that this is the largest value of $[Q]_1$. We begin by computing
 $$
      \frac{\partial [Q]_1}{\partial u_m} = -\big( u_m - \phibo(u_r) \big) - (u_m-u_l) < 0.
$$ 
Then we evaluate $[Q]_1$ when $u_m = p(u_l)$, where
$$
   p(u_l)  := u_l + \big(u_l-\phibSqr(u_l)\big) = u_l + \mysigma({\Npm}'), 
$$ 
and find
\begin{align*}
      [Q]_1(u_l,p(u_l), u_r) = &  \big( u_l-\phibSqr(u_l) \big) \, \big( \phibSqr(u_l) - \phibo(u_r) \big) \\
                     &  - \big( u_l-\phibSqr(u_l) \big) \, \big( 2 u_l - \phibSqr(u_l) - \phibo(u_r)\big)  \\
                 = & - 2 \big( u_l - \phibSqr(u_l) \big)^2 = -2 \mysigma({\Npm}')^2.       
\end{align*}
Therefore, for the unbounded set of states satisfying $u_m > p(u_l)$ 
(i.e. $\mysigma(\Rp)  > \mysigma(\Npm')$), we have $[Q]_1 < 0.$

For fixed $u_l$, the remaining set of states describing an RC-3 interaction form a bounded set 
in the $(u_m,u_r)$ plane defined by 
\begin{gather*}
   u_m <  p(u_l), \\
   \phis(u_m) < u_r < \phis(u_l).
\end{gather*}
Along the edge $u_m = p(u_l)$, $[Q]_1$ is negative and since $[Q]_1$ is decreasing
with respect to $u_m$, the largest value must occur along the edge $u_r = \phis(u_m)$.
We therefore compute the change in $[Q]_1$ along that edge, using $u_m$
as it's parameter.
\begin{align*}
  \frac{d}{d u_m} \Big( [Q]_1 &\big(u_l, u_m, \phis(u_m)\big)\Big) \\
  = &  \frac{d}{d u_m} \Big(
         \big( u_l-\phibSqr(u_l) \big) \, \big( \phibSqr(u_l) - \phibo \circ \phis(u_m) \big) \\
      &  \phantom{\frac{\partial}{\partial u_m} \Big( }           
            - (u_m - u_l) \, \big( u_m - \phibo \circ \phis(u_m) \big)
     \Big) \\
  = & \big( u_l-\phibSqr(u_l) \big) \, \Big( - \frac{d \phibo \circ \phis}{ d u} (u_m) \Big) \\
      &       -  \big( u_m - \phibo \circ \phis(u_m) \big) -(u_m - u_l) \, \Big( 1 - \frac{d \phibo \circ \phis}{ d u} (u_m) \Big)           
\end{align*}
Using the inequalities for the states and the fact that 
$\overline{\operatorname{Lip}}(\phibo \circ \phis) <  \overline{\operatorname{Lip}}(\phibSqr) < 1$, we
deduce that the terms above are negative and therefore, that $[Q]_1$ must attain it's
largest value for the smallest value of $u_m$, namely when $u_m = u_l = \phis(u_r)$
and $\sigma(\Rp) = 0$.

Assume now that $ \phis(u_l) >0$ and that $u_r$ could be positive but continue to fix $u_l$.
Then the set of states describing an RC-3 interaction form an infinite strip in the
$(u_m,u_r)$ plane defined by
$$
  u_l < u_m, \qquad \mbox{and} \qquad \phib(u_l) < u_r < \phis(u_l),
$$
where we have assumed the reasonable fact that $(\phis)' > 0$ at $u_r = \phis(u_l)$ (in fact, since
$(\phis)' < 0$ at the origin, we want to assume only one change of sign for $(\phis)'$ and
more general cases are similar). The same calculation as above shows that
$\partial [Q]_1 / \partial u_m < 0$ and therefore that the maximum value of
$[Q]_1$ must occur when $u_m = u_l$, i.e. $\sigma(R_{\text{in}})=0$.

\vskip.15cm 

\textbf{Case CR-4 :} ($\Cpm \Rm$)-(${\Npm}' {\Cm}'$).
The states satisfy \eqref{CR4_states}
$$
         \phibo(u_m) < \phibo(u_r)  < \phibSqr(u_l) < u_l,
$$
and therefore 
\begin{align*}
  [Q]_1 (u_l, u_m, u_r) = & \mysigma ( {\Npm}'  ) \mysigma({\Cm}' ) 
              - \mysigma( \Cpm)  \mysigma( \Rm ) \\
              = & \big(u_l - \phibSqr(u_l) \big) \big( \phibSqr(u_l) - \phibo(u_r) \big) \\
                  &  - \big(u_l-\phibo(u_m) \big) \big( \phibo(u_r) - \phibo(u_m) \big).
\end{align*}
If we let $u_r$ and $u_m$ approach each other while
maintaining the condition
$$
  u_r < \phis(u_l) < u_m,
$$
then $\phibo(u_r)-\phibo(u_m) \to 0$ while
$u_l - \phibo(u_m)$ remains bounded. Upon inspection,
it is clear that in this limit, $[Q]_1 \to  \mysigma ( {\Npm}'  ) \mysigma({\Cm}' )$.
We now demonstrate that this is the largest value of $[Q]_1$.
When $u_l$ is fixed, the admissible states belong to a rectangular domain
$$
      (u_m,u_r) \in \big[ \phis(u_l),0 \big] \times \big[ \phib(u_l), \phis(u_l)  \big].
$$
With respect to $u_m$ and $u_r$, we find
\begin{align*}
      \frac{\partial [Q]_1}{\partial u_m} =  \frac{d \phibo}{d u} (u_m) \, 
               \Big(  \big( \phibo(u_r) - \phibo(u_m) \big) + \big( u_l - \phibo(u_m)\big) \Big) < 0,
\end{align*}
and
\begin{align*}
      \frac{\partial [Q]_1}{\partial u_r} =
         - \frac{d \phibo}{d u} (u_r) \, 
               \Big(  \big( u_l - \phibSqr(u_l) \big) + \big( u_l - \phibo(u_m)\big) \Big) > 0.
\end{align*}
Taken together, these inequalities imply that the largest value of $[Q]_1$ occurs at
the corner of the domain where $u_m = u_r = \phis(u_l)$, i.e. when $\sigma(\Rm) = 0$.
On the other hand, the smallest value would occur along the boundary $u_r = \phib(u_l)$,
where one finds a negative quantity
\begin{align*}
  [Q]_1 =
  & (u_l - \phibSqr(u_l) ) (\phibSqr(u_l) - \phibSqr(u_l) ) 
  \\&                 - (u_l-\phibo(u_m)) (\phibSqr(u_l) - \phibo(u_m)).
\end{align*}

\textbf{Case CC-3 :} ($\Cp \CDown$)-(${\Npm}' {\CUp}'$).
As mentioned in the proof of Theorem \ref{thm:main}, the set of states are subdivided
into two subcases depending on the sign of $u_r$. When $u_r <0$, we identify
two possibilities \eqref{CC3_SUB1.1_states} and \eqref{CC3_SUB1.2_states}, namely 
\begin{align*} 
        \text{ either } \,  &   \phibo(u_r) < u_m < \phibSqr(u_l) < u_l,  
        \\
          \text{ or } \,  &  \phibo(u_r) <  \phibSqr(u_l) < u_m  < u_l, 
\end{align*}
where the second set of inequalities corresponds to a weak $\Cp$.
The change in $[Q]_1$ is then
\begin{align*}
  [Q]_1 (u_l, u_m, u_r) = & \mysigma ( {\Npm}'  ) \mysigma({\Cm}' ) 
              - \mysigma( \Cp)  \mysigma( \Cpm ) \\
              = & \big( u_l - \phibSqr(u_l) \big) \big(\phibSqr(u_l) - \phibo(u_r) \big) \\
                  & - (u_l-u_m) \big(u_m - \phibo(u_r) \big).
\end{align*}
Consider the function $B(\lambda) := \lambda (1-\lambda)$ whose maximum 
on the interval $[0,1]$ is attained at $\lambda=1/2$. There exists constants $C, \lambda$
and $\lambda'$ such that the change can be rewritten as
$$
   [Q]_1 = C \, \big( B(\lambda') - B(\lambda) \big). 
$$
It is therefore clear that $[Q]_1$ will be negative if and if $|\lambda - 1/2| < |\lambda'-1/2|$.
Unfortunately, this set cannot be described in a simple
manner in terms of the strengths of the waves.

\vskip.15cm 

\textbf{Case CN-3 :} ($\Cp \Npm$)-(${\Npm}' {\Cm}'$).
We identify two subcases
\begin{align*}
           \text{ either } \,  &  \phibSqr( u_m ) < u_m  < \phibSqr(u_l) < u_l, 
            \\
            \text{ or } \,  &  \phibSqr( u_m )   < \phibSqr(u_l) < u_m < u_l,
\end{align*}
where the second one occurs when the incoming shock $\Cp$ is
weak. The change can be written as
\begin{align*}
  [Q]_1 (u_l, u_m, u_r) = & \mysigma ( {\Npm}'  ) \mysigma({\Cm}' ) 
              - \mysigma( \Cp)  \mysigma( \Npm ) \\
              = & \big( u_l - \phibSqr(u_l) \big) \big( \phibSqr(u_l) - \phibSqr(u_m) \big) \\
                   & - (u_l-u_m) \big( u_m - \phibSqr(u_m) \big).
\end{align*}
Using a bit of algebra, we rewrite the change as
\begin{align}
  [Q]_1   = & \big( u_l - \phibSqr(u_l) \big) \big( \phibSqr(u_l) - \phibSqr(u_m) \big)   \nonumber \\
                      & - \big( u_l- \phibSqr(u_l) + \phibSqr(u_l) - u_m \big) \big( u_m - \phibSqr(u_m) \big) \nonumber \\
                   = & \big( u_l - \phibSqr(u_l) \big)  \, \big( \phibSqr(u_l) - u_m \big)            \nonumber \\
                       &  + \big( u_m - \phibSqr(u_l) \big) \, \big( u_m - \phibSqr(u_m) \big)   \nonumber  \\
                   = & - \big( u_m - \phibSqr(u_l) \big) \, \big( (u_l - u_m) 
                           - \big( \phibSqr(u_l) - \phibSqr(u_m) \big) \big),   \label{lemma:quad:RC3-main}
\end{align}
thus concluding that $[Q]_1< 0 $ as long as $\Cp$ is weak.
\end{flushleft}

\hfill $\Box$


\vskip.15cm 

{\it Proof of Theorem \ref{thm:quad}.}  
Throughout the proof, we use liberally the estimates derived in 
the proof of Theorem \ref{thm:main}. 

\vskip.15cm 

\begin{flushleft}
\textbf{Case RC-1 :} ($\Rp \CDown$)-(${\CDown}'$).
After  examining \eqref{RC1_SUB2_states}, it
is immediate that $\mysigma( \CDown) > \mysigma( {\CDown}' )$
and $W(\CDown) = W({\CDown}')$. Therefore, 
$$
    [Q]_1 = - \mysigma( \Rp) \mysigma (\CDown ) < 0
$$
and
\begin{align*}
   [Q]_2 & = W({\CDown}') \mysigma( {\CDown}') - W( \Rp ) \mysigma( \Rp )
                    - W(\CDown) \mysigma (\CDown ) \\
             & = - W(\Rp) \mysigma( \Rp) - 
                      W(\CDown) \big( \mysigma( \CDown) - \mysigma ({\CDown}') \big) < 0.       
\end{align*}

\textbf{Case RC-2 :} ($\Rp \Cpm$)-(${\Npm}' {\Rm}'$).
The conditions defining this case imply two subcases : 
\begin{align*}
     \text{either }  \, \phibSqr(u_l) < \phibo(u_r) < u_l < u_m, & 
    \\
   \text{or } \,   \phibSqr(u_l) < u_l < \phibo(u_r) < u_m. &
\end{align*}
Given that $\phis( u_m) < u_r < 0$, then property \eqref{contract} leads
to 
$$
0< \phibo(u_r) < \phibo \circ \phis (u_m) < \phibSqr(u_m) < u_m.
$$ 
This means that in the first subcase,
$$
\mysigma( {\Rm}') =  | \phibo(u_r) - \phibSqr(u_l)| < | \phibSqr (u_m) - \phibSqr(u_l) |
                                 <  | u_m - u_l | =    \mysigma( \Rp). 
$$ 
In that subcase, we also have
\begin{align*}
    \mysigma( \Cpm ) & = \mysigma( \Rp ) + | u_l - \phibo(u_r) |,
    \\
    \mysigma( {\Npm}' ) & = \mysigma( {\Rm}' ) + | u_l - \phibo(u_r) |, 
\end{align*} 
then $  \mysigma({\Npm}') < \mysigma( \Cpm) $. Similar arguments
prove the same inequalities in the second case. 

These inequalities therefore imply that
$$
    [Q]_1 = \mysigma({\Npm}') \mysigma( {\Rm}') - \mysigma(\Cpm) \mysigma(\Rp) < 0.
$$
Using the proposed definition
of {\sl weakly approaching waves},  we also deduce 
\begin{align*}
  [Q]_2  = &  W({\Npm}') \mysigma( {\Npm}' ) + W({\Rm}') \mysigma({\Rm}') 
                 - W(\Rp) \mysigma( \Rp) - W(\Cpm) \mysigma(\Cpm) \\
              = & 
                 - W({\Rm}') \big( \mysigma(\Rp) - \mysigma({\Rm}')  \big) 
                    - W({\Npm}') \, \big( \mysigma(\Cpm) - \mysigma( {\Npm}') \big) < 0.
\end{align*}

\textbf{Case RC-3 :} ($\Rp \CDown$)-(${\Npm}' {\CUp}'$).
Two subcases occur depending on the sign of $u_r$. 
If $u_r<0$, then \eqref{RC3_SUB1_states} holds and when $u_r >0$, 
then \eqref{RC3_SUB2_states} holds. In both of these cases
\be  \label{thm:quad:RC3-eq1}
  \mysigma( \CDown) - \mysigma(\Rp) = \mysigma({\CUp}') + \mysigma({\Npm}').
\ee
In Lemma \ref{lemma:quad}, we already showed that 
$[Q]_1$ could be positive.
For the other interaction term, we use \eqref{thm:quad:RC3-eq1} to verify 
\begin{align*}
[Q]_2 = & W( {\Npm}') \mysigma({\Npm}') + W({\CUp}') \mysigma( {\CUp}' ) 
            - W(\Rp) \mysigma(\Rp) - W(\CDown)\mysigma(\CDown) \\
          = & - W({\Npm}') \, \big(  \mysigma(\CDown) - \mysigma({\Npm}') - \mysigma({\CUp}') \big)
              - W(\Rp) \mysigma(\Rp) < 0.  
\end{align*}

\textbf{Case RN :} ($\Rp \Npm$)-(${\Npm}' {\Rm}'$).
The change $[Q]_1$ was studied in Lemma \ref{lemma:quad}.
For the second term, using the bounds on the wave strengths given in the
previous lemma, we again have a negative contribution
\begin{align*}
[Q]_2 = & W( {\Npm}') \mysigma({\Npm}') + W({\Rm}') \mysigma( {\Rm}' ) 
            - W(\Rp) \mysigma(\Rp) - W(\Npm)\mysigma(\Npm) \\
          = & - W({\Npm}) \, \big(   \mysigma(\Npm) - \mysigma({\Npm}')  \big)
              - W(\Rp) \, \big( \mysigma(\Rp) - \mysigma( {\Rm}')   \big) < 0.  
\end{align*}

\vskip.15cm 

\textbf{Case CR-1 :} ($\Cpm \Rm$)-(${\Cpm}'$).
In this case, the states satisfy $ \phibo(u_m) < \phibo(u_r) < \phibo( \phis(u_l) ) < u_l$ and
the wave strengths satisfy $\mysigma({\Cpm}') < \mysigma(\Cpm)$.  
Since only one wave is outgoing, $  [Q]_1 =   -  \mysigma(\Cpm) \mysigma( \Rm)  < 0$.  
On the other hand, it is easy to check that
\begin{align*}
[Q]_2 = & W( {\Cpm}') \mysigma({\Cpm}')  
            - W(\Cpm) \mysigma(\Cpm) - W(\Rm) \mysigma(\Rm) \\
          = & - W({\Cpm}) \, \big(   \mysigma(\Cpm) -  \mysigma({\Cpm}') \big)
              - W(\Rm)  \mysigma( {\Rm})  < 0.  
\end{align*}

\vskip.15cm 

\textbf{Case CR-2 :} ($\Cp \Rp$)-(${\Cp}'$).
This case is entirely classical so it is easy to check that 
\begin{align*}
  [Q]_1 = & - \mysigma(\Cp) \mysigma(\Rp) < 0, \\
  [Q]_2 = & - W(\Cp) \, \big(  \mysigma( \Cp)  - \mysigma({\Cp}')  \big)
                  - W(\Rp) \mysigma( \Rp ) < 0.
\end{align*}

\vskip.15cm 
\textbf{Case CR-3 :} ($\Cpm \Rm$)-(${\Npm}' {\Rm}'$).
A first subcase is defined by the additional condition $u_l < \phib(u_r)$ 
which provides   \eqref{CR3_SUB1_states} 
\be  
         \phibo(u_m) < \phibSqr(u_l) < u_l < \phibo(u_r).
\ee
The relative strengths of the waves are then seen to be
$   \mysigma( {\Rm}') < \mysigma( \Rm) $ and 
$  \mysigma( {\Npm}') < \mysigma( \Cpm)$.
In the second subcase, given by \eqref{CR3_SUB2_states}, these two
inequalities still hold. It is now easy to conclude
\begin{align*}
[Q]_1 = &  \mysigma({\Npm}') \mysigma({\Rm}')  -  \mysigma(\Cpm) \mysigma( \Rm) < 0, \\
[Q]_2 = & - W(\Cpm) \, \big( \mysigma( \Cpm ) - \mysigma( {\Npm}' )  \big)
                 - W(\Rm) \, \big( \mysigma(\Rm) - \mysigma( {\Rm}')  \big) < 0.
\end{align*}

\vskip.15cm 

\textbf{Case CR-4 :} ($\Cpm \Rm$)-(${\Npm}' {\Cm}'$).
The inequalities \eqref{CR4_states} defining this case suffice to show that
 $ \mysigma({\Npm}') + \mysigma( {\Cm}') \leq \mysigma( \Cpm) $.
In Lemma \ref{lemma:quad}, we showed that $[Q]_1 > 0$ but the second
term is negative 
$$
[Q]_2 = - W(\Cpm) \, \big( \mysigma(\Cpm) - \mysigma({\Npm}') - \mysigma( {\Cm}')  \big)
              - W(\Rm) \mysigma(\Rm) < 0.
$$

\vskip.15cm 

\textbf{Case CC-1 :} ($\Cp \CDown$)-(${\CDown}'$).
When $u_r$ is positive, all waves are classical and it is easy to show that 
$ [Q]_1 < 0$ and $[Q]_2 =0$. It is an exercise to see that when 
$u_r < 0$, then $[Q]_1<0$ is still negative and $[Q]_2$ vanishes.

\vskip.15cm 

\textbf{Case CC-2 :} ($\Cpm \CUp$)-(${\CDown}'$).
Two subcases appear depending on the sign of $u_r$ but each time,
$$
  [V] = \mysigma({\CDown}') - \mysigma(\Cpm) - \mysigma(\CUp) = 0.
$$
Only one wave is outgoing, so $[Q]_1<0$ and the previous identity implies
that  $[Q]_2 = 0$.

 \vskip.15cm 
\textbf{Case CC-3 :} ($\Cp \CDown$)-(${\Npm}' {\CUp}'$ ).
Since all waves involved are shocks and $[V]=0$, then
$[Q]_2=0$. 
\vskip.15cm 
\textbf{Case CN-1 :} ($\Cp \Npm$)-(${\Cpm}'$).
This is another simple case where the fact that $[V] =0$ and that 
$W(\cdot)$ is equal for all waves involved,  
implies  $[Q]_2 = 0$.

\vskip.15cm 
\textbf{Case CN-2 :} ($\Cpm \Nmp$)-(${\Cp}'$).
Same as Case CN-1.
\vskip.15cm 

\textbf{Case CN-3 :} ($\Cp \Npm$)-(${\Npm}' {\Cm}'$).
Two subcases occur but each time $[V] =0$. 
All waves are shocks so $[V]=0$ implies that $[Q]_2 = 0$.

\vskip.15cm 
\textbf{Case NC :} ($\Npm \CUp$)-(${\CDown}'$).
Same as Case CN-1.

\vskip.15cm 
\textbf{Case NN :} ($\Npm \Nmp$)-(${\Cp}'$).
Same as Case CN-1.

\end{flushleft}

\hfill $\Box$


\section{Global diminishing property for splitting-merging patterns}  
\label{sec:splitmerge}

In this section we show that, despite the fact that the 
quadratic interaction potential $Q_{\text{weak}}$ 
increases during interactions of type CR-3, RC-4, CN-3 and CC-3, this potential
is indeed {\sl strictly decreasing} 
globally in time for a large class of perturbations of crossing shocks. 
Hence, with the local bound we describe below, we provide the first steps towards an analysis of 
the global-in-time change of $Q_{\text{weak}}$ for arbitrary nonclassical entropy solutions. 
This section, therefore, provides a strong justification for the potential proposed in the previous section.

The {\sl splitting-merging solutions} considered now were introduced in LeFloch and Shearer \cite{LeFlochShearer},
where a modification of the total variation functional \cite{LeFloch-book} was shown to 
be strictly decreasing along the evolution of such splitting-merging solutions. 
The total variation functional $V$ presented in Section~\ref{sec:TV} also accomplishes this, 
but here we improve on those results by establishing a similar monotonicity 
result for the quadratic functional $Q_{\text{weak}}$.
Our analysis also brings to light some interesting aspects of splitting-merging solutions
that were not seen in \cite{LeFlochShearer}.

Splitting-merging solutions are, roughly speaking, perturbations of crossing shocks that
lead to the creation and destruction of a nonclassical shock. 
Such solutions contain two (classical and nonclassical) big waves that may merge together 
(as a classical shock) and also interact with (classical) small waves. 
A typical initial data for splitting-merging patterns 
is formed of 
\begin{enumerate}
\item[(i)] an isolated crossing shock with left- and right-hand states $u_-, u_+$ satisfying
$u_- > 0$ and $\phis( u_- ) < u_+$, but
$u_+ - \phis( u_- )$ small, 
\item[(ii)] followed, on the right-hand side, by a small rarefaction and a small shock. 
\end{enumerate}
The rarefaction is sufficiently strong that it has an interaction of type CR-4 
with the crossing shock, thereby leading to the creation of a pair of shock waves $\Npm, \CUp$. 
If the right-most shock is
sufficiently strong, then when it eventually interacts with $\CUp$ and the resulting shock will
begin to approach $\Npm$. The final interaction of type NC will involve 
$\Npm$ and the shock just described, thereby eliminating the nonclassical $\Npm$.
By adding more waves to the left and the right, this process of creation and destruction of
$\Npm$ can be repeated indefinitely.  

We consider a slightly more general configuration in the sense that we do not explicitly demand that 
a small shock on the right be responsible for the penultimate NC interaction. 
Fix some value $u^*> 0$ and define
\be 
 \label{defn:ubar}
    u_0^*(x) = \begin{cases}
                                        u^*,            & x < 0, \\
                                      \phis(u^*), & x >0.  
                         \end{cases}
\ee
Let $\theta_\eps$ be some function of locally bounded total variation and of oscillation bounded 
by some small positive $\eps$, i.e. 
$$
      \mysigma\big( \theta_\eps(x), 0 \big) < \eps, \qquad  x \in \RR.
$$
Without loss of generality, we may assume that $\theta_\eps$ is piecewise constant. 
Let $u^*_\eps$ be the nonclassical solution to the conservation law \eqref{conserv_law} with initial data
$u^*_0 + \theta_\eps$, as generated by the front-tracking method. 
Assuming the solution initially possesses a single isolated crossing shock located at the origin 
$x=0$, that is, assuming that 
${u^*_\eps(0-) > 0}$, 
we see that the crossing shock will be adjacent to many small classical 
shocks and rarefactions. After an interaction of type RC-3, CR-4, or CC-3, the small waves 
neighboring $\Cpm$ may lead to the creation of a pair of waves, a nonclassical shock $\Npm$
and a classical shock $\CUp$. After the creation of $\Npm$, 
the only types of interaction involving small waves incoming from the left of $\Npm$ are
RN and CN-3. The only types of interaction involving small waves and the shock $\CUp$, 
coming from either the left or the right,  are entirely classical (CC-1, RC-1 or CR-1). 
Moreover, no waves can cross $\CUp$ from the right or the left although the small waves
that reach $\Npm$ from the left, will cross and eventually reach $\CUp$. 
Therefore, the only way that the nonclassical shock $\Npm$ can be destroyed is 
if the shocks $\Npm$ and $\CUp$ change their speeds and eventually interact back together, 
leading us back to (a perturbation of) the original crossing shock.

These observations, in particular, 
imply that {\sl no waves can exit the domain} $\Omega$ bounded by the trajectories of $\Npm$ and $\CUp$. 
Our goal, in the present section, is to obtain a local bound on the change of the potential 
 $Q_{\text{weak}}$
relative only to the waves entering the domain $\Omega$. It should already be clear that the key here 
is comparing the total strength of the waves crossing $\Npm$ to the total strength of the waves 
terminating at $\CUp$.

Before stating our main result, we introduce some further notation. 
Let $t_0$ be a time of creation of a nonclassical wave $\Npm$ 
and denote by $t_1, t_2, \ldots, t_m$ the times of the next $m$ interactions between $\Npm$
and small waves $W_i$ on the right, and let $t_{m+1}$ be the time at which
$\Npm$ is destroyed from an interaction with the shock $\CUp$. Similarly,
let $\tilde{t}_i$ and $\overline{t}_i$ be the times at which an interaction occurs between
$\CUp$ and the left incoming waves $\widetilde{W}_i$ or
the right incoming waves $\overline{W}_i$, respectively. 
We define the {\sl total variation along the trajectory } $\Npm$ to be
\be
  \TV (\Npm) := \sum_{i=1}^m  \Big| \mysigma\big( \Npm(t_i+) \big)  - \mysigma\big( \Npm(t_i-) \big) \Big|,
\ee
and its {\sl signed variation} 
\be
  \SV (\Npm) :=  \mysigma\big( \Npm(t_{m+1}-) \big)  - \mysigma\big( \Npm(t_0+) \big).
\ee
Completely similar definitions also apply to the wave $\CUp$, but an additional decomposition 
exists by separating the contributions from the left- and the right-hand sides:
\begin{align*}
  \SV(\CUp) = & \bigg( \sum_{i=1}^{\widetilde{m}} \widetilde{s}_i \sigma ( \widetilde{W}_i )  \bigg)
                           + \bigg( \sum_{i=1}^{\overline{m}} \overline{s}_i \sigma ( \overline{W}_i ) \bigg) \\
                      =:    & \SV_L(\CUp) + \SV_R(\CUp),
\end{align*}
where $\widetilde{s}_i$ ($\overline{s}_i$) is $+1$ if $\widetilde{W}_i$ ($\overline{W}_i$) is a shock
and $-1$ otherwise.
For convenience, the strengths of the small wave $W_i$, before and 
after it has crossed $\Npm$ at some time $t_i$, are denoted by $W_i^-$ and $W_i^+$, respectively.

\vskip.15cm 

\begin{theorem}[Global diminishing property for splitting-merging patterns]
\label{thm:local_conserv}
Let $\phib$ be a kinetic function satisfying the properties (A1)--(A3), 
together with the mild requirement at some given reference state $u^*>0$.
$$
L^* :=   \big( \phibSqr )'(u^*) \in (1/2,1).
$$
Let $u^*_\eps$ be the nonclassical solution to the conservation law \eqref{conserv_law} with
initial data $u^*_0 + \theta_\eps$,
where $u^*_0$ is defined in \eqref{defn:ubar} and the perturbation $\theta_\eps$ is of locally bounded 
variation and of small amplitude, that is, $\| \theta_\eps \|_{L^\infty(\RR)} < \eps$.
Suppose that $u^*_\eps$ exhibits a splitting-merging pattern on the time interval 
$[t_0, t_{m+1}]$
along successive interactions with small waves (cf.~the notation above).
Provided $\eps$ is sufficiently small and the total effect of all waves on the classical 
shock $\CUp$ increases its total strength, that is,  
$$
\SV( \CUp) \geq 0,
$$
then the variation of the potential $Q=Q_{\text{weak}}$ within the domain $\Omega$ bounded by the trajectories of $\Npm$ and $\CUp$
is negative 
$$
[ Q]_1 \big|_{\Omega}<0.
$$ 
\end{theorem}

\vskip.15cm 

Recall that the expression $[ Q]_1$ was defined earlier in Section~4. 
In the above theorem, only contributions within the domain 
$\Omega$ are considered. 
Two important remarks about our assumptions are in order: 
\begin{enumerate}

\item
  The condition on $L^*$ is satisfied for a wide range of fluxes, entropies and nonlinear 
  diffusion-dispersion models, see Theorems~2.3 and 3.3 in \cite{LeFloch-book}.
  In particular, it is true for the cubic flux $f(u) = u^3-au$, the entropy $U(u) = u^2/2$
  and all ratios of diffusion over dispersion. The non-trivial aspect of this condition is the
  lower bound $ L^*>1/2$ since the upper bound is always satisfied.

\item Our main assumption $\SV(\CUp) \geq 0$ requires that the total effect of the
interaction of all waves on $\CUp$ is to increase its strength.
In fact, this is always the case for the perturbations of splitting-merging solutions
within the setting \cite{LeFlochShearer}. In our slightly more general setting though,
 waves crossing through $\Npm$
will change the critical state $\phis(u_l)$ and therefore could conceivably lead to an
NC interaction even if $\CUp$ interacted only with rarefactions on the right ($\SV(\CUp) < 0$).
Intuitively, one would like to show that the existence of an NC interaction at time $t_{m+1}$
implies that the wave $C$ became strong enough to change course
and therefore that $\SV(\CUp) \geq 0$ is satisfied. This is the subject of ongoing research.

\end{enumerate}

  \
  
It is interesting to check that the condition $\SV( \CUp) \geq 0$ of Theorem \ref{thm:local_conserv}
 is in fact necessary, by studying a simple situation involving
 a single rarefaction $R$ crossing the nonclassical shock $\Npm$,
  and one classical shock $C'$ interacting
 from the right with $\CUp$ (assuming that $\phis$ is monotonically decreasing, say). 
 Suppose that the strength of the rarefaction wave after crossing 
 $\Npm$ is $\nu_R$ and that the change in $\Npm$ during this interaction
 is $\nu_L$. In Lemma \ref{lem:local_conserv} below, we will check that $\nu_L < \nu_R$
 and that $\nu_L+\nu_R$ is the strength of the original rarefaction.
 So, assuming that $C'$ interacts with $\CUp$ before the rarefaction does, then the total change in
 $Q$ is given by :
\be
\nonumber
\aligned
 & [Q]_1  |_{\Omega} \, = \, && + \sigma(\Npm) \sigma(\CUp)             &&&   \text{ (due to CR-4),} 
 \\
 &                                      && - \nu_L \big( \sigma(\Npm) + \nu_R)   &&&  \text{ (due to RN),}   
 \\
 &                                      && -  \sigma(\CUp) \sigma(C')            &&& \text{ (due to CC-1),}  
 \\
 &                                      && -  \big( \sigma( \CUp ) + \sigma(C') \big) \nu_R   
                  &&& \text{ (due to RC-1),} 
 \\
 &                                      && -  \big( \sigma(\Npm) - \nu_L   \big) \, 
                                              \big( \sigma(\CUp) + \sigma(C') -\nu_R \big)        
                                                &&& \text{ (due to   NC).}
\endaligned
\ee 
The last term can be rewritten as
$$
     -  \sigma(\Npm) \sigma(\CUp) -\sigma(\Npm) \big( \sigma(C') -\nu_R \big)
     + \nu_L  \big( \sigma(\CUp) +    \sigma(C')  \big)-  \nu_L  \nu_R.
$$
The contribution from the RC-1 interaction can control the third term.
Unfortunately, if $\nu_R$ is large with respect to $\sigma(C')$, i.e. $\SV( C') < 0$, then the
CC-1 term cannot be used to control the second term.


For the proof of Theorem~\ref{thm:local_conserv} we will need the following two lemmas. 
The interest of the first lemma is to make more precise the (mainly linear) dependence of the change $[Q]_1$ 
in terms of the incoming wave $\sigma(W_i^-)$. The second lemma provides an estimate which
closely relates the signed variations of $\CUp$ and $\Npm$. 

 \vskip.15cm 
 
 \begin{lemma}[Interactions with the nonclassical shock]  
 \label{lem:local_conserv}
 Consider interactions RN and CN-3 at the time $t_i$ involving a left-incoming weak
 wave $W_i$ and the nonclassical shock $\Npm$. Then, there exists a positive constant
 $L_i < 1$ such that
$$
\aligned 
\sigma( W_i^+ ) 
& = L_i \sigma( W_i^-), 
\\
\sigma\big( \Npm( t_i+) \big) 
& =  \sigma\big( \Npm( t_i-) \big) + s_i (1- L_i )\sigma( W_i^-), 
\\
[Q]_1                 
& = - (1-L_i) \sigma(W_i^-) \, \Big( \sigma\big( \Npm(t_i -) \big)- s_i L_i \sigma(W_i^-)  \Big) < 0,
\endaligned
$$
where $s_i$ is $+1$ if $W_i$ is a shock and $-1$ otherwise. Moreover, one has 
\be  
\label{eq3:lem_local}
      L_i = L^* + \BigO{ \eps }.  
\ee
\end{lemma}

\vskip.15cm 

{\it Proof.} We consider only a RC-3 interaction,
since the calculations for an incoming rarefaction are similar and have been
essentially treated in Lemma \ref{lemma:quad}, equation \eqref{lemma:quad:RC3-main}.  
When $W_i$ is small, then the states are ordered as in \eqref{CN3_SUB2_states}, namely
$$
 \phibSqr(u_m^i) < \phibSqr(u_l^i) < u_m^i < u_l^i. 
 $$
Then, we have 
$$ 
     \sigma( W_i^+ )  = \Big| \phibSqr (u_l^i) - \phibSqr (u_m^i) \Big| =: L_i | u_l^i - u_m^i | 
                                    = L_i \sigma( W_i^-),
$$ 
with $L_i < 1$ since $\phibSqr$ is a contraction. 
On the other hand, relation \eqref{eq3:lem_local} is obvious since 
the $u_m^i, u_l^i$ belong to a neighborhood of $u^*$ of size $\epsilon$. 

Finally, the outgoing nonclassical shock has strength
\begin{align*}
   \sigma\big( \Npm(t_i+) \big) & = \big| u_l^i - \phibSqr (u_l^i) \big|  \\
                & = |u_l^i - u_m^i| + \big| u_m^i - \phibSqr (u_m^i) \big| - \Big| \phibSqr (u_l^i) - \phibSqr (u_m^i) \Big| \\
                & = \sigma\big( \Npm( t_i-) \big) +  (1- L_i )\sigma( W_i^-),
\end{align*}
while the change $[Q]_1$ takes the form 
\begin{align*}
  [Q]_1 & = \sigma \big(  \Npm(t_i+) \big) \sigma( W_i^+) -  \sigma \big( \Npm(t_i-) \big) \sigma( W_i^- )  \\
             & = \Big(  \sigma \big( \Npm( t_i-) \big) +  (1- L_i )\sigma( W_i^-) \Big) \, L_i \sigma(W_i^-) 
                 - \sigma \big( \Npm(t_i-) \big) \sigma( W_i^- ) \\
             & = - (1-L_i) \,  \sigma(W_i^-)  \, \Big(  \sigma\big( \Npm( t_i-) \big) - L_i \sigma( W_i^-) \Big) < 0.       
\end{align*}
\hfill $\Box$


\vskip.15cm 

\begin{lemma} [Property of the signed variations]
  \label{lem:local_SV}
Setting $\lambda^* := (1-L^*)/L^*$, one has  
\be
 \Big|  \lambda^* \SV_L(\CUp) - \SV(\Npm) \Big| \leq \BigO{\eps} \TV(\Npm).
\ee
\end{lemma}

\vskip.15cm 

{\it Proof.}  
Using Lemma \ref{lem:local_conserv}, we compute
\begin{align*}
\SV(\Npm) = & \sum_{i=1}^m s_i (1-L_i) \sigma(W_i^-) \\
          = & \sum_{i=1}^m s_i (1-L^*) \sigma(W_i^-)  + \sum_{i=1}^m s_i (L^*-L_i) \sigma(W_i^-)  
                       \\
           = & \lambda^* \sum_{i=1}^m s_i L^* \sigma(W_i^-)  + \sum_{i=1}^m s_i (L^*-L_i) \sigma(W_i^-),
\end{align*}  
thus
$$
\SV(\Npm) = \lambda^* \SV_L(\CUp) + (1+\lambda^*) \sum_{i=1}^m s_i (L^*-L_i) \sigma(W_i^-). 
$$
For each index $i$, $L^* - L_i = \mathcal{O}_i( \eps )$, so if we take $\BigO{\eps} = \max_i  \mathcal{O}_i( \eps ) $
and use $\underline{L} := \min_i L_i$, then the difference 
is   
\begin{align*}
    \Big|  \lambda^* \SV_L(\CUp) - \SV(\Npm) \Big| & \leq 2 \BigO{\eps} \sum_{i=1}^m \big| \sigma(W_i^-) \big| \\
     & \leq \frac{2}{\underline{L}} \BigO{\eps} \TV(\Npm) = \BigO{\eps} \TV(\Npm). 
\end{align*}
\hfill $\Box$

\ 

{\it Proof of Theorem \ref{thm:local_conserv}.}  To help the reader understand the wave interactions (and cancellations) 
in this proof, we begin with a few preliminary comments. 
Simply put, the final NC interaction should provide a quadratic term 
which is, later, cancelled by a similar quadratic term when the two waves merge. One may naively expect 
the sum 
$$
- \sigma(\Npm (t_{m+1}-)) \, \sigma( \CUp (t_{m+1}-))
+ \sigma(\Npm (t_0+)) \, \sigma( \CUp (t_0+))
$$
to be negative. Of course, the cumulative strength of the changes during the 
interactions with the small waves  must also be taken into account. 
The difference between the strength of the initial and final waves $\Npm, \CUp$ is measured by
the signed variation $\SV$ along those two shocks. 
Our proof below shows that, along the trajectories $\Npm$ and $\CUp$, the
change $[Q]_1$ is negative and proportional to the total variation $\TV(\Npm) + \TV(\CUp)$. 
The total variation being larger than the signed variation, after further analysis 
one can conclude that $[Q]_1 |_{\Omega} < 0$.

The key technical information is provided by Lemma~\ref{lem:local_SV},
which implies that, up to a quantity of order $\BigO{\eps} \TV(\Npm)$, 
we have $\SV(\Npm) = \lambda^* \, \SV_L(\CUp)$ and, in particular, that the signed variations have the same sign.
We now have all the tools necessary to proceed with the proof of Theorem~\ref{thm:local_conserv}.

The perturbation $\theta_\eps$ has bounded oscillation and therefore it can only alter the right-hand state of
$\CUp$ by an amount $\eps$. The small waves entering $\Omega$ through $\Npm$
only alter its strength by $| (1-L_i) \, \sigma( W_i^- ) | < \eps$. 
In both cases, we expect that everywhere along their trajectories,
\begin{align*}
    \sigma \big( \Npm( t ) \big) & > \sigma \big( \Npm( t_0+)\big) - 2\eps, \\
    \sigma \big( \CUp( t ) \big) & > \sigma \big( \CUp( t_0+)\big) - 2\eps.
\end{align*}

Ignoring the negative contribution to $[Q]_1$ coming from  the interaction
that generated $\Npm$ (which, anyway, can be arbitrarily small), and neglecting also 
all classical interactions inside $\Omega$ (for which $[Q]_1 \leq 0$ and possibly  $0$), 
we have
\begin{align*}
[Q]_1 \Big|_{\Omega} \leq &  + \sigma\big(\Npm (t_0+) \big) \, \sigma \big( \CUp (t_0+) \big) \\
                            &-  \sum_{i=1}^m (1-L_i) \sigma(W_i^-) \, \Big( \sigma\big( \Npm(t_i -) \big)- s_i L_i \sigma(W_i^-)  \Big) \\
                            & - \sum_{i=1}^{\widetilde{m}} \sigma ( \widetilde{W}_i ) \, \sigma \big( \CUp (\widetilde{t}_i-) \big) 
                             - \sum_{i=1}^{\overline{m}}  \sigma ( \overline{W}_i ) \, \sigma \big( \CUp (\overline{t}_i-) \big) \\
                             & -  \sigma \big(\Npm (t_{m+1}-) \big) \, \sigma \big( \CUp (t_{m+1}-) \big), 
\end{align*}                          
and after simplification
\begin{align*}
[Q]_1 \Big|_{\Omega}                           
                           \leq & + \sigma \big(\Npm (t_0+) \big) \, \sigma \big( \CUp (t_0+) \big) \\
                           & - \Big( \sigma\big( \Npm(t_0+) \big) - 3 \eps \Big) \, \TV( \Npm) 
                            - \Big( \sigma\big(  \CUp(t_0+) \big) - 2 \eps \Big) \, \TV( \CUp) \\
                            & - \Big( \sigma\big( \Npm (t_0+) \big) + \SV(\Npm)  \Big) \, 
                                      \Big( \sigma\big(\CUp (t_0+) \big) + \SV( \CUp) \Big).
\end{align*}

Introduce constants $\eta_N, \eta_C \in [0,1]$ to be determined 
later in the proof. The terms can then be split in the following way.
\begin{align}  \label{eq:fund_decomp}
[Q]_1 \Big|_{\Omega} \leq &  - \Big( (1-\eta_N) \sigma\big( \Npm(t_0+) \big) - 3 \eps \Big) \, \TV( \Npm)  \\
                            & - \Big( (1-\eta_C) \sigma\big( \CUp(t_0+) \big) - 2 \eps \Big) \, \TV( \CUp)  \nonumber  \\
                            & -  \sigma\big(\Npm (t_0+)\big) \, \Big( \eta_N \TV(\Npm) + \SV(\CUp)  \Big)    \nonumber   \\
                            &   -   \sigma\big(\CUp (t_0+)\big) \, \Big( \eta_C \TV(\CUp) + \SV( \Npm) \Big) 
                              - \SV(\Npm) \, \SV( \CUp).      \nonumber 
\end{align}
If $\eta_N, \eta_C$ are known a priori, the first two terms are negative (and can be neglected) 
by taking $\eps$  sufficiently small. Clearly, the main difficulty now lies in the sign of
$\SV(\Npm) $ and $\SV( \CUp)$. If both are positive, then all the terms are negative and we have
proved our result. Recall that $\SV( \CUp) \geq 0$, by assumption, so that it suffices to consider 
the (only possibly unfavorable) case $\SV(\Npm) < 0$. 

Lemma \ref{lem:local_SV} implies that 
we can assume that $\SV_L(\CUp)$ is also negative, provided 
we accept a small correction term $\BigO{\eps} \TV(\Npm)$. 
From this, we deduce that $\SV_R(\CUp) > 0$ and in fact, necessarily
\begin{equation*} 
          \big| \SV_L(\CUp)  \big|  < \big| \SV_R(\CUp) \big|. 
\end{equation*}
We immediately note that any correction terms of order $\BigO{\eps} \TV(\Npm)$ can be included
into the first term of the decomposition \eqref{eq:fund_decomp} and, therefore, 
taking a smaller $\eps$ if necessary, we can make the new term negative:
$$
    - \Big( (1-\eta_N) \sigma\big( \Npm(t_0+) \big) - \BigO{\eps} \Big) \, \TV( \Npm) < 0.
$$
This fact is used below without further comment.

Using Lemma \ref{lem:local_SV} we have $| \SV(\Npm)| \approx \lambda^* | \SV_L(\CUp)| \leq \lambda^* \TV(\CUp)$, 
and the fourth term in \eqref{eq:fund_decomp} can be written as
\begin{align*}
  -  \sigma  (\CUp (t_0+)) \, \Big( &\eta_C \TV(\CUp)  + \SV( \Npm) \Big) \\
   \leq & -   \sigma(\CUp (t_0+)) \, \Big( \eta_C \TV(\CUp) - \lambda^* \, \TV( \CUp) \Big)
            + \BigO{\eps} \TV(\Npm). 
\end{align*}
Taking $\eta_C = \lambda^*$ suffices to guarantee that the term vanishes.
We now re-organize the third and fifth terms of \eqref{eq:fund_decomp} 
to take advantage of the known signs of the signed variations, as follows: 
\begin{align}  
\nonumber 
A := &
-  \sigma \big(\Npm (t_0+)\big) \, \Big(  \eta_N \TV(\Npm) + \SV(\CUp)  \Big)     - \SV(\Npm) \, \SV( \CUp) \\
       = &    -  \sigma\big(\Npm (t_0+) \big) \, \Big( \eta_N \TV(\Npm) +\frac{\eta_N}{2} \SV_L(\CUp)  \Big)   \nonumber\\
          &    -   \sigma\big(\Npm (t_0+)\big) \, \Big( \big(1- \frac{\eta_N}{2} \big) \SV_L(\CUp) + \SV_R(\CUp)  \Big) \nonumber \\
          &      - \SV(\Npm) \, \SV_L( \CUp)  - \SV(\Npm) \, \SV_R( \CUp) + \BigO{\eps} \TV(\Npm).  \nonumber 
\end{align}
Observe that we have the following lower bounds
$$
     \big(1-\frac{\eta_N}{2} \big)  \SV_L(\CUp) + \SV_R(\CUp) \geq \frac{\eta_N}{2} \Big|  \SV_L(\CUp) \Big|
$$
and upper bounds
$$
    \Big|  \SV_L(\CUp) \Big| = \frac{1}{\lambda^*} \Big|  \SV(\Npm) \Big| +  \BigO{\eps} \TV(\Npm)
                                               \leq  2  \TV(\Npm)  +  \BigO{\eps} \TV(\Npm).
$$
Using the two previous  bounds, we find 
\begin{align}   \nonumber  
A
       \leq &    -  \sigma\big(\Npm (t_0+) \big) \, \Big( \eta_N \TV(\Npm) - \eta_N \TV(\Npm)  \Big)   
                -   \sigma\big(\Npm (t_0+)\big) \, \frac{\eta_N}{2} \Big|  \SV_L(\CUp) \Big|  \nonumber \\
             &      - \SV(\Npm) \, \SV_L( \CUp)  - \lambda^* \SV_L (\CUp) \,  \SV_R( \CUp) + \BigO{\eps} \TV(\Npm) \nonumber \\
    \leq   &    -   \Big|  \SV_L(\CUp) \Big| \, \Big(  - \lambda^* \,  \SV_R( \CUp) 
                                                                                    + \frac{\eta_N}{2}   \sigma\big(\Npm (t_0+)\big) \Big)  \nonumber \\
             &   - \SV(\Npm) \, \SV_L( \CUp) + \BigO{\eps} \TV(\Npm). \nonumber 
\end{align}
For any fixed value of  $\eta_N \in (0,1)$, the first term is negative because $|\SV_R(\CUp)| < \epsilon$
can be assumed a priori small. The second term is negative and the last one can be incorporated into the first term
of \eqref{eq:fund_decomp}. This completes the proof of Theorem \ref{thm:local_conserv}.
\hfill $\Box$


\section{Quadratic interaction potential (part 2)}
\label{sec:quad2}

We continue our investigation of interaction functionals adapted to nonclassical shocks. 
The proposal now is to also involve a weight-function in front of the quadratic terms
which is chosen to be proportional to the difference in (normalized) wave speeds. Such a weighting was 
used successfully by Iguchi and LeFloch \cite{IguchiLeFloch} as well as by Liu and Yang \cite{LiuYang}
in their analysis of 
general hyperbolic systems of conservation laws; 
see also Liu \cite{Liu} for earlier pioneering work on systems.   

The interaction functional proposed now for nonclassical entropy solutions shares the same monotonicity properties as 
the functional $Q_{\text{weak}}$,
in the sense that it fails for precisely the same interactions. In preparation for the
study of systems of conservation laws, a complete study of the scalar case requires an
examination of the following functional. 

So, we introduce 
$$
Q_{\text{pos}}(u(t)) := \sum_{x<y} \theta(x,y) \sigma(x) \sigma(y), 
$$
where $\sigma(x), \sigma(y)$ are the generalized wave strengths defined earlier, 
$$
\sigma(x) 
=
\begin{cases}
\sigma(S), & S \text{ is a wave located at } x, 
\\
0,         & \text{ no wave located at } x. 
\end{cases}
$$
The potential now contains a weight function $\theta(x,y)$ defined  by 
$$
\theta(x,y) := 
\begin{cases}
C_* \, \big( \ah(x) - \ah(y) \big)^+,    & \text{same monotonicity,}
\\
1, & \text{otherwise,} 
\end{cases}
$$
where the constant $C_*>0$ will be chosen to be sufficiently small and 
$(k)^+ := \max\{ k,0\}$. 
More precisely, here we say that two waves have the same monotonicity if 
after mapping all states to the positive region the waves are either both increasing or both decreasing. 
The definition of $\ah$ is as follows: 
$$
\ah(x) := {f(\uh_x^+) - f(\uh_x^-) \over \uh_x^+ - \uh_x^-}, 
$$
where 
$$
\uh_x^+ := \begin{cases}
u_x^+,          & u_x^+ \geq 0, 
\\
\phibo(u_x^+),          & u_x^+ <0. 
\end{cases}
$$
This definition allows us to enforce that comparisons only occur between state variables
belonging to the same region of convexity, i.e. between positive values.  
Note that, in the present formalism, rarefaction are possibly interacting with each other. 

\vskip.15cm

\begin{theorem}[``Positive interaction'' potential for nonclassical shocks]
\label{thm:positive}
Let $\phib$ be a kinetic function satisfying the properties (A1)--(A3). 
Consider the weighted interaction functional $Q_{\text{pos}}$ evaluated
for a front-tracking approximate solution $u$ to \eqref{conserv_law}, whose
discontinuities satisfy \eqref{defn_phib}, and with initial data 
$u_0 \in \operatorname{BV}(\mathbb{R})$. If  
$C_*$ is taken sufficiently small, then $Q_{\text{pos}}$ is strictly decreasing 
during all interactions except in the cases RC-3, CR-4, CC-3 and CN-3
which involve a small wave interacting with a crossing shock.
\end{theorem}

\vskip.15cm 

{\it Proof.} We distinguish between several cases of interactions. For the functional 
$Q=Q_{\text{pos}}$ we introduce a new form of the decomposition \eqref{decompo},
$$
       \big[ Q \big(u(t)\big) \big] =: [Q]_1 + \sum_y  \sigma(y) \, B(y),
$$
where the sum over y excludes waves involved in the interaction at time $t$.
Given that the waves not involved in the interaction are potentially arbitrary,
$[Q]$ will be negative if and only if $[Q]_1$ and $B(y)$ are negative. As usual, we will
rely heavily on the inequalities described in the proofs of Theorems \ref{thm:main} and
\ref{thm:quad} as well as of Lemma \ref{lemma:quad}.

\vskip.15cm 

\begin{flushleft}
{\bf Case RC-1:}  ($\Rp \CDown$)-(${\CDown}'$).
In this case, a single wave is outgoing and we have $[Q]_1 < 0$. 
In the proof of Theorem \ref{thm:main}, we observed that  the states of interest satisfy 
\begin{align*}
       \text{ either } \,  & u_r < u_l < u_m, 
      \\
          \text{ or } \,  &  \phibo(u_r) < u_l < u_m.  
\end{align*}
For both sets of inequalities, $\sigma(\CDown) < \sigma(\CDown')$ holds 
and the normalized speeds are ordered
$$
\ah(\CDown') < \ah(\CDown) < \ah (\Rp).
$$ 
Suppose there is a shock at some point $y$ with speed $\Lambda$, then 
$$
B(y) = C_* \, \big(\ah(\CDown') - \Lambda \big)^+ \sigma(\CDown') 
            -  C_* \,  \big(\ah(\CDown) - \Lambda \big)^+ \sigma(\CDown) - \sigma(\Rp), 
$$
which by the monotonicity just observed must be negative. On the other hand, if 
a rarefaction is located at $y$, then the estimate on the wave strengths shows that 
$$
B(y) = \sigma(\CDown') 
            -  \sigma(\CDown) - C_* \,  \big(\ah(\Rp) - \Lambda \big)^+ \sigma(\Rp) < 0. 
$$

\vskip.15cm 

\textbf{Case RC-2:} ($\Rp \Cpm$)-(${\Npm}' {\Rm}'$). Our study of this interaction in 
Theorem \ref{thm:quad} provides  
$$
\sigma(\Npm') < \sigma(\Cpm), \qquad \sigma(\Rm') < \sigma(\Rp). 
$$
In turn, it follows that 
$$
[Q]_1 = \sigma(\Npm') \, \sigma(\Rm') - \sigma(\Cpm) \, \sigma(\Rp) < 0. 
$$
Next, to control $B(y)$ we first consider the relevant speeds.
When $\phibo(u_r) < u_l$, then inequalities \eqref{RC2_SUB1_states} are valid and
$$
\ah(\Rm') < \ah(\Npm) < \ah (\Cpm) < \ah(\Rp),
$$
while when $\phibo(u_r) > u_l$ we have 
$$
\ah(\Nmp') < \ah(\Rm') < \ah (\Rp) < \ah(\Cpm). 
$$
Now, suppose that a shock wave is located at some point $y$, then 
$$
\aligned
B(y) = 
& C_* \,  \big(\ah(\Npm') - \Lambda \big)^+ \sigma(\Npm') + \sigma(\Rm') 
\\
&            -  C_* \,  \big(\ah(\Cpm) - \Lambda \big)^+ \sigma(\Cpm) - \sigma(\Rp), 
\endaligned
$$
which is negative in view of the above inequalities on the speeds and strengths. 
A similar argument applies if a rarefaction is located at $y$. 

\vskip.15cm 

\textbf{Case RC-3 :} ($\Rp \CDown$)-(${\Npm}' {\CUp}'$).
This case is more delicate, due to the fact that a nonclassical shock is generated from classical waves. 
The shock strengths satisfy the relation \eqref{thm:quad:RC3-eq1} which implies 
\be   \label{thm:+:RC3_eq1}   
      \sigma (\Npm') + \sigma(\CDown') < \sigma(\CDown). 
\ee
We distinguish between two cases \eqref{RC3_SUB1_states} and \eqref{RC3_SUB2_states} 
depending on the sign of $u_r$. An inspection of the normalized wave speeds associated to the 
states in \eqref{RC3_SUB1_states}  and \eqref{RC3_SUB2_states}  demonstrates that
$$
     \ah(\CUp') < \ah(\Npm') < \ah(\Rp) \qquad \mbox{and} \qquad
     \ah(\CUp') < \ah( \CDown) < \ah(\Rp).
$$

We analyzed an interaction term in Lemma \ref{lemma:quad} similar to
\begin{align*}
[Q]_1 
 = C_*  \big( \ah(\Npm') - \ah(\CUp') \big)^+ \, \sigma({\Npm}') \, \sigma(\CUp') 
    - \sigma(\Rp) \sigma(\CDown),
\end{align*}
and we saw that it was possible to take $\sigma(\Rp)$ vanishingly small, while
maintaining a non-vanishing product $\sigma({\Npm}') \sigma(\CUp')$. Therefore, in general
we cannot claim that $[Q]_1$ is negative.

 We now show that in the limit as $\sigma(\Rp) \to 0$, $B(y)$ is again positive
for some speeds $\Lambda$. Assuming a shock located at $y$,  the general form of $B$ is  
$$
\aligned
B(y) = 
& C_* \,  \big(\ah(\Npm') - \Lambda \big)^+ \, \sigma(\Npm') 
       + C_* \,  \big(\ah(\CUp') - \Lambda \big)^+ \, \sigma(\CUp') 
\\
&       - C_* \,  \big(\ah(\Cp) - \Lambda \big)^+ \, \sigma(\Cp) - \sigma(\Rp). 
\endaligned
$$
If $\phis(u_l) < 0$ and inequalities \eqref{RC3_SUB1_states} hold, then take
$u_r \to \phis(u_l)$ and $u_m \to u_l$. In this limit, $ \ah(\CUp') < \ah(\CDown) < \ah(\Npm') $
and the quantities $\ah(\CUp')$ and $\sigma(\Npm')$ are functions only of $u_l$. Therefore
taking $\Lambda \in [ \ah(\CDown), \ah(\Npm')]$, we find that for those limiting values
$$
    B(y) = C_* \,  \big(\ah(\Npm') - \Lambda \big)^+ \, \sigma(\Npm') > 0.
$$

Now, if a rarefaction is located at $y$, then \eqref{thm:+:RC3_eq1} suffices to show
$$
B(y) =   \sigma(\Npm') + \sigma(\CUp')     -  \sigma(\CDown) 
                   - C_* \,  \big(\ah(\Rp) - \Lambda \big)^+ \, \sigma(\Rp) < 0.
$$

\vskip.15cm 

\textbf{Case RN :} ($\Rp \Npm$)-(${\Npm}' {\Rm}'$). 
As we saw in the proof of Theorem \ref{thm:main}, there are two subcases defined by
either 
\begin{align*}
       \text{ either } \,  &  \phibo \circ \phib (u_l) < \phibo \circ \phib (u_m) < u_l < u_m, 
        \\
       \text{ or } \,  &  \phibo \circ \phib (u_l) < u_l < \phibo \circ \phib (u_m) < u_m. 
\end{align*}
The normalized speeds satisfy respectively In Case A the speeds satisfy 
\begin{align*}
       \ah(\Rm') < \ah(\Npm') < \ah(\Npm) < \ah(\Rp), &  \mbox{ or} \\
       \ah(\Npm') < \ah(\Rm') < \ah(\Rp) < \ah(\Npm).
\end{align*}
In both subcases, we have
$\sigma(\Npm') < \sigma(\Npm)$ and $\sigma(\Rm') < \sigma(\Rp)$, 
and the inequalities on the speeds $\ah(\Npm')  < \ah(\Npm) $ and
$\ah({\Rm}')  < \ah(\Rp)$.

Using the observations above, we verify that
$$
[Q]_1 = \sigma(\Npm') \sigma(\Rm') - \sigma(\Rp) \, \sigma(\Npm) < 0. 
$$ 
If a shock is located at $y$, then the term
$$
\aligned
B(y) = 
& C_* \,  \big(\ah(\Npm') - \Lambda \big)^+ \, \sigma(\Npm') + \sigma(\Rm') 
+ C_* \,  \big(\ah(\Npm) - \Lambda \big)^+ \, \sigma(\Npm) - \sigma(\Rp)
\endaligned
$$ 
is negative. When a rarefaction is located at $y$, then 
$$
\aligned
B(y) = 
& \sigma(\Npm') + C_* \,  \big( \ah(\Rm') - \Lambda \big)^+ \, \sigma(\Rm') 
  - \sigma(\Npm) - C_* \,  \big( \ah(\Rm) - \Lambda \big)^+ \, \sigma(\Rm) 
\endaligned
$$ 
is again negative.

\vskip.15cm 

\textbf{Case CR-1 :} ($\Cpm \Rm$)-(${\Cpm}'$).
Only one wave is outgoing, so we trivially have $[Q]_1 < 0$. 
The states of interest satisfy 
$$
        \phibo(u_m) < \phibo(u_r) < u_l, 
$$
and therefore the normalized speeds are ordered 
$$
\ah(\Rm) < \ah(\Cpm) < \ah(\Cpm').  
$$
Recall also that $\sigma(\Cpm') < \sigma(\Cpm)$. Assume now that a shock is located at $y$,
then 
$$
\aligned
B(y) = 
& C_* \, \big( \ah(\Cpm') - \Lambda \big)^+ \, \sigma(\Cpm') 
       - C_* \, \big( \ah(\Cpm) - \Lambda \big)^+ \, \sigma(\Cpm) 
         - \sigma(\Rm). 
\endaligned
$$
The least favorable case corresponds to the situation where $\Lambda = \ah(\Cpm)$, 
since only one negative term remains when $\Lambda > \ah({\Cpm}')$ and
$B(y)$ is increasing with respect to $\Lambda$ when $\Lambda < \ah(\Cpm)$. The important
term in $B(y)$ when $\Lambda = \ah(\Cpm)$ is the difference in speeds
\begin{align*}
        \big( \ah(\Cpm') - \ah(\Cpm) \big)^+  
             & = \bigg| \frac{f(u_l) - f(\uh_r)}{u_l-\uh_r} -   \frac{f(u_l) - f(\uh_m)}{u_l-\uh_m} \bigg| \\
          & \leq C \, \big| \uh_r - \uh_m \big| = C \sigma(\Rm).
\end{align*}
The constant $C = C(u_l)$ is bounded above since the total variation
of the solution $u$ remains bounded. If we therefore assume the a priori bound
\be  \label{thm:+:condition2}
       C_* C \TV(u) < 1,
\ee
then 
$$
\aligned
B(y) \leq 
& C_* C \, \sigma(\Rm) \, \sigma(\Cpm') - \sigma(\Rm)
= -  \sigma(\Rm)  \big( 1 - C_* C \,\, \sigma(\Cpm') \big) < 0.
\endaligned
$$
To complete the proof, we assume that a rarefaction is now located at $y$, and immediately see that
$$
B(y) = \sigma(\Cpm') - \sigma(\Cpm) - C_* \, (\ah(\Rm) - \Lambda)^+ \, \sigma(\Rm) < 0. 
$$

\vskip.15cm

\textbf{Case CR-2 :} ($\Cp \Rp$)-(${\Cp}'$). 
The waves are entirely classical and the calculations are identical to previous Case CR-1.
In particular, we again have to impose condition \eqref{thm:+:condition2} on $C_*$. 

\vskip.15cm 

\textbf{Case CR-3 :} ($\Cpm \Rm$)-(${\Npm}' {\Rm}'$).
In the proof of Theorem \ref{thm:main}, we identified two subcases distinguished by the 
inequalities \eqref{CR3_SUB1_states} and \eqref{CR3_SUB2_states} :
\begin{align*}
  \text{ either } \,  & \phibo(u_m) < \phibo \circ \phib(u_l) < u_l < \phibo(u_r), 
  \\ 
      \text{ or } \,  &  \phibo(u_m) < \phibo \circ \phib(u_l) <  \phibo(u_r) < u_l. 
\end{align*}
In both subcases we have 
$$
\sigma(\Npm') < \sigma(\Cpm), \qquad \sigma(\Rm') < \sigma(\Rm). 
$$
In the first subcase, the normalized speeds satisfy 
$$ 
\ah(\Cpm) < \min\big( \ah(\Rm), \ah(\Npm') \big) < \max\big( \ah(\Rm), \ah(\Npm') \big)  < \ah(\Rm'), 
$$
while in the second, 
$$ 
\ah(\Rm) < \min\big( \ah(\Cpm), \ah(\Rm') \big) < \max\big( \ah(\Cpm), \ah(\Rm') \big) < \ah(\Npm').
$$
The key point is that in both subcases, $\ah(\Cpm) < \ah({\Npm}')$ and
$\ah(\Rm) < \ah({\Rm}')$.

Using the inequalities satisfied by only the strengths, we obtain 
$$
[Q]_1 = \sigma(\Npm') \sigma(\Rm') - \sigma(\Cpm) \sigma(\Rm) < 0.
$$
When a shock is located at $y$, we can verify that 
$$
B(y) = C_* \, \big( \ah(\Npm') - \Lambda \big)^+ \, \sigma(\Npm') + \sigma(\Rm')
     - C_* \, \big(\ah(\Cpm) - \Lambda \big)^+ \, \sigma(\Cpm) - \sigma(\Rm).
$$  
The worst case occurs when $\Lambda = \ah(\Cpm) < \ah(\Npm')$, and then we find
that there exists a constant $C$ such that
\begin{align*}
   \ah({\Npm}') - \ah(\Cpm) & = \bigg| \frac{f(u_l) - f\big( \phibSqr(u_l) \big)}{u_l - \phibSqr(u_l)}  
                                     -   \frac{f(u_l) - f(\uh_m)}{u_l-\uh_m} \bigg|  \\
                                     & \leq  C \big(  \phibSqr(u_l) - \uh_m \big)
                                            =   C \big( \sigma(\Rm) - \sigma(\Rm') \big).
\end{align*}
We again impose a condition of the form \eqref{thm:+:condition2} on $C_*$ and find
$$
\aligned
B(y) = 
& C_* \, \big( \ah(\Npm') - \ah(\Cpm) \big)^+ \, \sigma(\Npm') - \big( \sigma(\Rm) - \sigma(\Rm') \big) 
\\
= 
& - \big( 1 - CC_* \, \sigma(\Npm') \big) \, \big( \sigma(\Rm) - \sigma(\Rm') \big) < 0.
\endaligned
$$

Suppose next that there is a rarefaction at $y$, then 
$$
\aligned
B(y) = 
& \sigma(\Npm')  + C_* \, \big( \ah(\Rm') - \Lambda \big)^+ \, \sigma(\Rm') 
         -  \sigma(\Cpm) - C_* \, \big( \ah(\Rm) - \Lambda \big)^+ \, \sigma(\Rm).
\endaligned
$$ 
The largest value of $B(y)$ occurs when $\Lambda = \ah(\Rm)$ and then the coefficient
of the $\sigma(\Rm')$ term satisfies 
\begin{align*}
  \ah(\Rm') - \ah(\Rm) & = \bigg| \frac{f(\uh_r) - f\big( \phibSqr(u_l) \big)}{\uh_r - \phibSqr(u_l)}  
                                     -   \frac{f(\uh_r) - f(\uh_m)}{\uh_r-\uh_m} \bigg|  \\
                                     & \leq  C \big(  \phibSqr(u_l) - \uh_m \big) 
                                = C \, \big( \sigma(\Cpm) - \sigma(\Npm') \big).
\end{align*}
The same arguments as before apply and again allow us to conclude that $B(y) <0$.

\vskip.15cm

\textbf{Case CR-4 :} ($\Cpm \Rm$)-(${\Npm}' {\Cm}'$).
Here, we have 
$$
     \phibo(u_m) < \phibo(u_r) < \phibo \circ \phibo(u_l) < u_l.
$$
Now, the normalized speed satisfy 
$$
\ah(\Rm) < \ah(\Cm') < \ah(\Npm'), 
\qquad 
\ah(\Rm) < \ah(\Cpm) < \ah(\Npm')
$$
and for the strengths we have 
$$
\sigma(\Npm') + \sigma(\Cm') < \sigma(\Cpm).
$$
For the interaction, we obtain a term
$$
[Q]_1 = C_* \, \big( \ah(\Npm') - \ah(\Cm') \big)^+ \, \sigma(\Npm') \sigma(\Cm') - \sigma(\Cpm) \sigma(\Rm)
$$
similar to the one analyzed in Lemma \ref{lemma:quad}. There we saw that it
was possible to let $\sigma(\Rm) \to 0$ while keeping both 
$\sigma(\Npm')$ and $\sigma(\Cm)$ bounded and non-vanishing. Therefore,
this term is positive for a large class of interactions.

When a shock is located at $y$, the general term to look at is
\begin{align*}
   B(y) = & C_* \big( \ah(\Npm') - \Lambda \big)^+ \sigma(\Npm') + C_* \big(\ah(\Cm') - \Lambda\big)^+ \sigma(\Cm') \\
         & - C_* \big( \ah(\Cpm)-\Lambda\big)^+ \sigma(\Cpm) - \sigma(\Rm).
\end{align*}
Taking again the limit as $\sigma(\Rm) \to 0$, we find 
$  \ah(\Cm') < \ah(\Cpm) < \ah(\Npm') $. Moreover, when $\Lambda = \ah(\Cpm)$
we obtain expressions $\sigma(\Npm')$ and $\ah(\Npm') - \ah(\Cpm)$ that are non-zero and functions
only of $u_l$. Therefore
$$
    B(y) = C_* \big( \ah(\Npm') - \Lambda \big)^+ \sigma(\Npm') > 0.
$$

When a rarefaction is located at $y$, then
$$
  B(y) = \sigma(\Npm') + \sigma(\Cm') - \sigma(\Cpm) 
          - C_*\big( \ah(\Rm) - \Lambda \big)^+ \sigma(\Rm) < 0.
$$

\vskip.15cm 

\textbf{Case CC-1 :} ($\Cp \CDown$)-(${\CDown}'$).
There are two subcases determined by the sign of $u_r$ :
$$
  \phibo(u_r) < u_m < u_l, \qquad \mbox{or} \qquad u_r < u_m < u_l.
$$
 In both subcases, we have the identities
\begin{gather*}
      \sigma(\CDown') =  \sigma( \Cp) + \sigma( \CDown), \\
      \ah( \CDown) < \ah(\CDown') < \ah(\Cp).
\end{gather*}
Only one wave is outgoing so $[Q]_1<0$.

Suppose a shock is located at $y$, then the largest value of
$$
   B(y) = C_* \big( \ah(\CDown') - \Lambda \big)^+ \sigma(\CDown') 
           - C_* \big( \ah(\Cp) - \Lambda \big)^+ \sigma( \Cp) 
           - C_* \big(  \ah(\CDown) - \Lambda \big)^+ \sigma(\CDown)
$$
occurs in the interval $\Lambda \in[ \ah(\CDown), \ah(\CDown')]$, leaving us with
$$
   B(y) = C_* \big( \ah(\CDown') - \Lambda \big)^+ \sigma(\CDown') 
           - C_* \big( \ah(\Cp) - \Lambda \big)^+ \sigma( \Cp).
$$
In the subcase where $u_r >0$, the definition of $\sigma$ and $\ah$ allows us
to find
\begin{align*}
   B(y) = & C_* \, \Big(  \big( f(u_l) - f(u_r) \big) - \Lambda (u_l - u_r) 
                    -  \big( f(u_l) - f(u_m) \big) + \Lambda (u_l - u_m)  \Big)
                      \\
            = & C_* \big(  \ah(\CDown) - \Lambda \big)  \, \sigma(\CDown) \leq 0.       
\end{align*}
If $u_r < 0$, then similar computations show that $B(y)$ is still negative.

When a rarefaction is located at $y$, then
$$
    B(y) = \sigma(\CDown') - \sigma(\Cp) - \sigma(\CDown) = 0.
$$

\vskip.15cm 

\textbf{Case CC-2 :} ($\Cpm \CUp$)-(${\CDown}'$).
Two subcases appear depending on the sign of $u_r$. We will assume throughout that
$u_r >0$, leaving us with the states
$$
  u_r < \phibo(u_m) < u_l,
$$
and therefore the two identities
\begin{gather*}
      \sigma(\CDown') =  \sigma( \Cpm) + \sigma( \CUp), \\
      \ah( \CUp) < \ah(\CDown') < \ah(\Cpm).
\end{gather*}
We recognize immediately the same identities as in the Case CC-1. As quick
verification shows that the same arguments demonstrate that $[Q]_1<0$ and
$[Q]_2<0.$

\vskip.15cm 

\textbf{Case CC-3 :} ($\Cp \CDown$)-(${\Npm}' {\CUp}'$ ).
We can  distinguish between two subcases :
\begin{align*}
       \text{ either } \,  &  \phibo(u_r) < u_m < \phibSqr(u_l) < u_l, 
     \\
       \text{ or } \,  &  \phibo(u_r) < \phibSqr(u_l)  < u_m  < u_l.  
\end{align*}
In both subcases, one obtains
\be   \label{thm:+:CC3_eq2}
    \ah(\CDown) < \ah(\Npm') \qquad \mbox{and} \qquad \ah(\CUp') < \ah(\Cp).
\ee

Unfortunately, the conditions above on the shock speeds do not preclude the possibility that
$\ah(\Cp) < \ah(\Npm')$  which would imply that
\begin{align*}
    B(y) = & C_* \big( \ah(\Npm') - \Lambda \big)^+ \sigma(\Npm') 
                   + C_* \big( \ah(\CUp') - \Lambda \big)^+ \sigma(\CUp') \\
                & - C_* \big( \ah(\Cp) - \Lambda \big)^+ \sigma(\Cp) 
                     - C_* \big( \ah(\CDown) - \Lambda \big)^+ \sigma(\CDown) 
\end{align*}
was positive. On the other hand, because $[V]=0$ we have that 
$B(y)=0$ when a rarefaction is located at $y$.

The sign of $[Q]_1$ is also less than helpful. Using the ideas seen 
in Lemma \ref{lemma:quad} for the same term, one notices that
when $u_m$ equals the transitional state $ \phibSqr(u_l)$ between the two subcases,
then the change
$$
     [Q]_1 = C_* \big( \ah(\Npm') - \ah(\CUp') \big)^+ \sigma(\Npm') \sigma(\CUp')
                   - C_* \big( \ah(\Cp) - \ah(\CDown)  \big)^+  \sigma(\Cp) \sigma(\CDown)
$$
vanishes because the quadratic terms and the coefficients
are equal. It is clear that for some set
of waves close to this transitional value, $[Q]_1$ will also be positive.

\vskip.15cm 

\textbf{Case CN-1 :} ($\Cp \Npm$)-(${\Cpm}'$).
The states are $   \phibSqr(u_m) < u_m < u_l $
and these provide us with the basic inequalities
\begin{gather*}
      \sigma(\Cpm') =  \sigma( \Cp) + \sigma( \Npm), \\
      \ah( \Npm) < \ah(\Cpm') < \ah(\Cp).
\end{gather*}  
Under these conditions, the analysis is similar to the CC-1 case.

\vskip.15cm 

\textbf{Case CN-2 :} ($\Cpm \Nmp$)-(${\Cp}'$).
The important states are ordered $\phib(u_m) < \phibo(u_m) < u_l$
and therefore
\begin{gather*}
      \sigma({\Cp}') =  \sigma( \Cpm) + \sigma( \Nmp), \\
      \ah( \Nmp) < \ah({\Cp}') < \ah(\Cpm).
\end{gather*}  
Again, the analysis is similar to the CC-1 case.

\vskip.15cm

\textbf{Case CN-3 :} ($\Cp \Npm$)-(${\Npm}' {\Cm}'$).
There are two subcases defined by the ordered set of states
\begin{align*}
   \text{ either } \,  & \phibSqr(u_m) < u_m < \phibSqr(u_l) < u_l, 
    \\
   \text{ or } \,     & \phibSqr(u_m) < \phibSqr(u_l)  < u_m < u_l,
\end{align*}
where the second case corresponds to a weak incoming shock $\Cp$.
It is easy
to verify that the shock strengths verify 
\be  \label{thm:+:CN3_eq1}
    \sigma(\Npm) < \sigma(\Npm') \qquad \mbox{and} \qquad \sigma(\CUp') < \sigma(\Cpm).
\ee
In both subcases, one obtains the unfortunate inequalities
\be   \label{thm:+:CN3_eq2}
    \ah(\Npm) < \ah(\Npm') \qquad \mbox{and} \qquad \ah(\CUp') < \ah(\Cpm).
\ee

The analogue to 
$$
   [Q]_1 = C_* \big( \ah(\Npm') - \ah(\Cm') \big)^+ \sigma(\Npm') \sigma(\Cm')
                - C_* \big( \ah(\Cp) - \ah(\Npm)  \big)^+  \sigma(\Cp) \sigma(\Npm)
$$
for $Q_{\text{weak}}$ was analyzed in Lemma \ref{lemma:quad} and seen to be 
positive for sufficiently strong waves $\Cp$. For the functional
$Q_{\text{pos}}$, $[Q]_1=0$ when $u_m = \phibSqr(u_l)$ because at that
critical value, all inequalities in \eqref{thm:+:CN3_eq1} and \eqref{thm:+:CN3_eq2}
are equalities. Without going into the details, it easy to verify that this implies that
for a non-empty set of waves, $[Q]_1$ will be positive.

Suppose a shock is located at $y$, then
\begin{align*}
    B(y) = & C_* \big( \ah(\Npm') - \Lambda \big)^+ \sigma(\Npm') 
                   + C_* \big( \ah(\Cm') - \Lambda \big)^+ \sigma(\Cm') \\
                & - C_* \big( \ah(\Cp) - \Lambda \big)^+ \sigma(\Cp) 
                     - C_* \big( \ah(\Npm) - \Lambda \big)^+ \sigma(\Npm). 
\end{align*}
When we are in the first subcase (a strong $\Cp$), the change
$B(y)$ is clearly positive if $\Lambda \in [\ah(\Cp), \ah(\Npm')]$.

Finally, suppose a rarefaction is located at $y$, then
$$
 B(y) = \sigma(\Npm') + \sigma(\Cm') - \sigma(\Cp) - \sigma(\Npm) = 0.
$$

\vskip.15cm

\textbf{Case NC :} ($\Npm \CUp$)-(${\CDown}'$).
The defining states are $u_r < \phibSqr(u_l) < u_l$ and we have
\begin{gather*}
      \sigma(\CDown') =  \sigma( \Npm) + \sigma( \CUp), \\
      \ah( \CUp ) < \ah(\CDown') < \ah(\Npm).
\end{gather*}  
See the CC-1 case.

\vskip.15cm

\textbf{Case NN :} ($\Npm \Nmp$)-(${\Cp}'$).
The analysis of the NC case also applies here since NN interactions
are limiting situations obtained as $ u_r \to \phib(u_m)$.  
 
\end{flushleft}

\hfill $\Box$ 


\section*{Acknowledgments}
The first author (M.L.) was partially supported by the Natural Sciences and Engineering Research Council (NSERC)
of Canada. The second author (PLF) was partially supported by the Centre
National de la Recherche Scientifique (CNRS) and the Agence
Nationale de la Recherche (ANR) through the grant 06-2-134423
entitled {\sl ``Mathematical Methods in General Relativity''} (MATH-GR).


\bibliographystyle{plain} 

\begin{thebibliography}{1}

\bibitem{AbeyaratneKnowles}  \auth{R. Abeyaratne and J.K. Knowles,}
Kinetic relations and the propagation of phase boundaries in solids,
Arch. Rational Mech. Anal. 114 (1991), 119--154.

\bibitem{LeFlochBedjaoui}  \auth{N. Bedjaoui and P.G. LeFloch,}
Diffusive-dispersive travelling waves and kinetic relations V. Singular
diffusion and nonlinear dispersion,
Proc. Royal Soc. Edinburgh 134A (2004), 815-844.

\bibitem{IguchiLeFloch}  \auth{T. Iguchi and P.G. LeFloch,}
Existence theory for hyperbolic systems of conservation laws with general flux-functions, 
Arch. Rational Mech. Anal. 168 (2003), 165--244. 

\bibitem{Glimm} \auth{J. Glimm,} 
Solutions in the large for nonlinear hyperbolic systems of equations,
Comm. Pure Appl. Math. 18 (1965), 697--715.
 
\bibitem{LeFloch}  \auth{P.G. LeFloch,}
Propagating phase boundaries: formulation of the problem and existence via the Glimm scheme,
Arch. Ration. Mech. Anal. 123 (1993), 153--197.

\bibitem{LeFloch-book}  \auth{P.G. LeFloch,}
{\sl Hyperbolic systems of conservation laws:
The theory of classical and nonclassical shock waves},
Lectures in Mathematics, ETH Z\"urich, Birkh\"auser, 2002.

\bibitem{Dafermos72}  \auth{C.M. Dafermos,}
Polygonal approximation of solutions to the initial-value problem for a conservation law,
J. Math. Anal. Appl. 38 (1972), 33--41.

\bibitem{LeFlochMohammadian} \auth{P.G. LeFloch and M. Mohammadian,} 
Why many shock wave theories are necessary.  Fourth-order models, kinetic functions, and equivalent equations, 
J. Comput. Phys. 227 (2008), 4162--4189. 

\bibitem{LeFlochShearer}  \auth{P.G. LeFloch and M. Shearer,}
Nonclassical Riemann solvers with nucleation,
Proc. Royal Soc. Edinburgh 134A (2004), 941--964.

\bibitem{Liu} \auth{T.-P. Liu,} 
{\sl Admissible solutions of hyperbolic conservation laws,} 
Memoir Amer. Math. Soc. 30, 1981.

\bibitem{LiuYang} \auth{T.-P. Liu and T. Yang,}
Weak solutions of general systems of hyperbolic conservation laws, 
Commun. Math. Phys. 230 (2002), 289--327.

\bibitem{Slemrod}  \auth{M. Slemrod,}
Admissibility criteria for propagating phase boundaries in a van der Waals fluid,
Arch. Rational Mech. Anal. 81 (1983), 301--315.

\bibitem{Truskinovsky} \auth{L. Truskinovsky,}
Kinks versus shocks,
in ``Shock induced transitions and phase structures in general media'',
R. Fosdick, E. Dunn, and M. Slemrod ed., IMA Vol. Math. Appl., Vol.~52,
Springer-Verlag, New York (1993), pp.~185--229.



\end{thebibliography}

\newcommand{\auth}{\textsc}

\end{document}